\documentclass[twocolumn]{autart}

\usepackage{graphicx}
\usepackage{mathtools}
\usepackage{todonotes}
\usepackage{txfonts}
\usepackage{fontawesome }
\usepackage{derivative}
\usepackage{booktabs}
\usepackage{macros}
\usepackage{empheq}
\usepackage{lscape}
 \usepackage{xcolor,colortbl}
\definecolor{lavender}{rgb}{0.9, 0.9, 0.98}
\usepackage{adjustbox}
\usepackage{caption}
\usepackage{subcaption}
\usepackage{tikz}
\usepackage{pgfplots}  
\usepackage{enumitem}
\setitemize{nolistsep}
\setenumerate{nolistsep}
\usepackage{utfsym}
\usepackage{url}
\usepackage{bm}

\usepackage{tabularray}       
\UseTblrLibrary{booktabs}     
\usepackage{adjustbox}        
\usepackage{caption}          
\usepackage{float}
\usepackage{multirow}
\usepackage{algorithm}
\usepackage{algpseudocode} 

\usepackage{url}

\allowdisplaybreaks

\begin{document}
\setlength{\abovedisplayskip}{2pt}
\setlength{\belowdisplayskip}{2pt}
\setlength{\abovedisplayshortskip}{2pt}
\setlength{\belowdisplayshortskip}{2pt}

\begin{frontmatter}

 \title{Algorithmic detection of false data injection attacks in cyber-physical systems\thanksref{footnoteinfo}}

\thanks[footnoteinfo]{Souvik Das was supported by a PMRF grant from the Ministry of Human Resource Development, Govt. of India, during his tenure in IIT Bombay, India. Debasish Chatterjee acknowledges partial support of the SERB MATRICS grant MTR/2022/000656.}

\author[SD]{Souvik Das}\ead{souvik.das@angstrom.uu.se},
\author[AG]{Avishek Ghosh}\ead{avishek\_ghosh@iitb.ac.in},
\author[AG]{Debasish Chatterjee}\ead{dchatter@iitb.ac.in}

\address[SD]{Division of Signals and Systems, Department of Electrical Engineering, Uppsala University, Sweden.}
\address[AG]{Center for Systems \& Control Engineering, Indian Institute of Technology Bombay, Powai, Mumbai - 400076, India}



\begin{abstract}
This article introduces an anomaly detection based algorithm (AD-CPS) to detect \emph{false data injection attacks} that fall under the category of data deception/integrity attacks, but \emph{with arbitrary information structure}, in cyber-physical systems (CPSs) modeled as stochastic linear time-invariant systems. The core idea of this data-driven algorithm is based on the fact that an honest state (one not compromised by adversaries) generated by the CPS should concentrate near its weighted empirical mean of the immediate past samples. As the first theoretical result, we provide non-asymptotic guarantees on the \emph{false positive error} incurred by the algorithm for attacks that are \emph{\(2\)-step honest}, referring to adversaries that act intermittently rather than successively. Moreover, we establish that for adversaries possessing a certain minimum energy, the \emph{false negative error} incurred by AD-CPS is low. Extensive experiments were conducted on partially observed stochastic LTI systems to demonstrate these properties and to quantitatively compare AD-CPS with an optimal CUSUM-based test.
 \end{abstract}


\end{frontmatter}

\section{Introduction}
\label{sec:intro}
Growing autonomy and interconnectivity in cyber-physical systems (CPSs) have enabled them to deliver better performance, efficiency, and enhanced fault resilience in many practical applications. However, these advancements come at a cost, resulting in abundant exposure to adversarial attacks. Some attacks target the cyber layer of CPSs consisting of the cryptographic and communication layers, and once breached leave the entire physical layer of the CPS open to manipulation. Securing the physical layer of a CPS is, therefore, imperative.
Such attacks on the \emph{physical layer} of CPSs may deceive the honest entities in the CPS, devise plans to tamper with the nominal operations of CPSs, corrupt data on communication channels, and simultaneously, may even compromise the entire system \cite{ref:DK-13,ref:EK-14}. Consequently, the automatic and fast online detection of anomalies caused by these adversaries is crucial; failure to detect anomalies may have severe consequences, especially when these underlying attacked CPSs are critical to our daily livelihood. An emergent need to develop sophisticated detection schemes to detect the presence of these adversaries, therefore, exists. 
See \cite{ref:RM-IRC-14,ref:HH-JY:16,ref:JG-DU-AC-etal-18,ref:DD-GLH-YX-etal-18,ref:HS-VG-KHJ:22} for surveys of such techniques and their applications in various contexts.

It is a widely accepted trend among researchers that an adversary attacking a CPS should be resource-constrained \cite{ref:AT-IS-HS-15,ref:KP-AT-MC-PP-17} and must be described by an \emph{attack model} \cite{ref:AT-DP-HS-KHJ-12}. While it is reasonable to assume such an attack model for adversaries up to a certain extent, an intelligent adversary, on the other hand, may play malicious actions arbitrarily, not per the belief of the designer: they may collude or choose to wait for a long time to gather information and learn the underlying system's parameters before launching a full-blown attack based on their accumulated data. 
Therefore, detecting the presence of such \emph{near-omnipotent} adversaries is crucial, difficult, and hardly studied in the literature. 

Some recent results developed in \cite{ref:SD-PD-DC-23} attempt to bridge that gap by partially addressing these questions for large classes of history-dependent and \emph{intelligent actuator attacks}. Two key features of the results of \cite{ref:SD-PD-DC-23} involve \((a)\) \emph{model free adversaries}, and \((b)\) \emph{almost sure detectability of the presence of adversaries}; however, the established results are asymptotic and difficult to verify. In contrast, our focus here is on algorithms with non-asymptotic theoretical guarantees capable of rapid detection of these adversaries. 
We consider adversaries capable of attacking sensors or hijacking communication channels by injecting additive disturbances at intermittent time instances; see \S\ref{subsec:attack model} for more details. Moreover, \emph{we do not impose restrictions on these adversarial disturbances in the sense that they can either be a sequence of deterministic vectors, functions of the prior history of the system, or vectors randomized according to some unknown attack/policy adopted by the adversary.}
This class of algorithms falls under the ambit of anomaly detection; see \cite{ref:VC-AB-VK-09,ref:DS-21} for a survey of such techniques; and they may also be a good candidate for rapid fault detection (both adversarial and non-adversarial) in CPSs due to its data-driven nature. In that context, \cite{ref:AW-HJ-76} is perhaps one of the earliest works where tools from statistics were extensively used, laying the groundwork for subsequent research in this field.
A good detection algorithm not only needs to detect anomalies caused by adversaries but also should not trigger when the CPS is operating normally. Two types of errors are relevant --- \emph{False positive error (FPE) or Type-\(\text{I}\) error} when the detection algorithm asserts the presence of adversaries even though they are not present or active in reality, and \emph{False negative error (FNE) or Type-\(\text{II}\) error} when the algorithm fails to detect the presence of adversaries even though they are present --- and we assess the algorithm's performance in terms of these errors.
\subsection*{Related work}
Let us discuss several recent works on sequential detection schemes and anomaly-detection-based methods.
%
Asymptotic guarantees on the false alarm and detection rates were established in \cite{ref:YM-RC-BS-13}. 
Finite-sample guarantees for analyzing the false positive error over a time horizon were established in \cite{ref:PH-MP-RV-AA-20} in the context of sensor-switching attacks; however, no comments were made regarding false negative errors. The article \cite{ref:DL-SM-20} presented a detection scheme employing the Wasserstein metric and provided ``high-confidence'' performance guarantees over finite samples, while \cite{ref:DU-JR-HS-22} supplied finite-sample guarantees on the false alarm rate. Along similar lines, an anomaly detection test based on the Kantorovich distance was established in \cite{ref:MM-RKK-22}, although without theoretical guarantees. For systems with multiplicative noise, the article \cite{ref:VR-BJG-JR-THS-22} developed a residual-based detection test, but no theoretical guarantees were presented. An optimal CUSUM-based detection test with a watermarking scheme was developed in \cite{ref:AN-AT-AA-SD-23} and asymptotic performance guarantees were reported, but the test is tailored to a Gaussian noise model and a specific choice of the attack model. More recently, \cite{ref:SD-AG-DC-24} established a sequential detection test based on the separability of state trajectories and \cite{ref:SD-PD-DC-23} provided guarantees on false positive errors; however, that detection scheme required continuous monitoring of the system's nominal behavior.

\subsection*{Contributions of the article}
\begin{itemize}[leftmargin=*, label = \(\circ\)]
        \item \textbf{An anomaly detection based framework:} This article presents a \emph{general framework} for a family of anomaly identification based detection algorithms that employ a simple threshold-based test at the \emph{signal level} to detect a broad range of false data injection attacks \cite{ref:YM-EG-AC-BS-10,ref:BS-PRK-16,ref:PH-MP-RV-AA-17,ref:CZB-FP-VG-17,ref:TYZ-DY-20,ref:MG-KG-WL-20,ref:SP-ATK-22,shinohara2025data,shinohara2025detection,tosun2025kullback,feng2025false}, and \cite{ye2025decentralized,anand2025analysis} on sensor channels in stochastic LTI systems subject to general \emph{sub-Gaussian noise}. In particular, we present a detection algorithm --- AD-CPS --- that operates on data generated by a completely observable system. At its core, AD-CPS is conceptually motivated by \cite{ref:AS-BN-VYS-ZS:22} and relies on the premise that an \emph{honest state vector} (states that are not corrupted) should concentrate near the empirical weighted mean of the recent past states. 
    
        \item \textbf{Key features of AD-CPS:} Our detection algorithm is compatible with sub-Gaussian noise and \emph{does not} require complete knowledge of the process noise distribution. This leads to an inherent \emph{distributional robustness}, in contrast to most recent works \cite{chen2025moving,xiao2025guaranteed,li2024secure,jiang2025attack} that rely on Gaussian data. Moreover, AD-CPS \emph{does not} assume the underlying attack model a priori nor the knowledge of the attack policy employed by these adversaries.

        \item \textbf{Guarantee on the false positive error:} The algorithm is equipped with \emph{provable ``with high probability'' guarantees} regarding the false alarm rate under certain sufficient conditions on the class of adversaries, as detailed in Theorem \ref{th:error_incurred}. Unlike works based on sequential tests such as CUSUM and its variants or \(\chi^2\)-based failure tests, we provide a closed-form expression for the threshold that ensures a low false-positive error. Moreover, Theorem \ref{th:error_incurred} does not depend on, nor distinguishes between, the specific attack policies employed by adversaries --- a feature that also reveals an inherent trade-off between system stability and algorithm performance, which is studied extensively through numerical experiments.


        \item \textbf{Guarantee on false negative error:} The algorithm is also supported by \emph{``with high probability'' guarantees} on the false negative error (mis-detection rate). Establishing such guarantees is generally challenging and has been only sparsely discussed in the literature.

        \item \textbf{Detailed experimentations:} Extensive experiments were conducted on partially observed systems, which are more general and realistic than fully observable ones. Key aspects of the algorithm, such as threshold selection and the trade-off between false positive and false negative errors, have been thoroughly investigated and reported here. These studies underscore the generality and robustness of our approach. Comparisons with an optimal CUSUM test \cite{ref:AN-AT-AA-SD-23} have demonstrated comparable performance.

\end{itemize}



\subsection*{Organization}
The article unfolds as follows:
\S\ref{sec:ls} formally states the question that we ask in this article. Our algorithm --- AD-CPS --- is presented in \S\ref{sec:mrg}, followed by a series of remarks highlighting several technical aspects of the algorithm. \S\ref{sec:erros} and \S\ref{sec:fn_err} consist of the main results of this article --- Theorem \ref{th:error_incurred} and Theorem \ref{th:t2_error_incurred}. The proofs of Theorem \ref{th:error_incurred} and Theorem \ref{th:t2_error_incurred} are relegated to Appendix \ref{appen:t1_error} and Appendix \ref{appen:t2_error}, respectively. All the auxiliary results needed to prove Theorem \ref{th:error_incurred} and Theorem \ref{th:t2_error_incurred} are provided in Appendix \ref{appen:proofs}.

\subsection*{Notations and symbols}
Let \((\ps,\tsigalg,\PP)\) be a (sufficiently rich) probability space, and we assume that all random elements are defined on \((\ps,\tsigalg,\PP)\). The realization of a \(d\mbox{-}\)dimensional random vector \(g\) at \(\omega \in \ps\) defined on \((\ps,\tsigalg,\PP)\) is given by \(g(\omega)\). In the rest of the article, we omit these details for brevity, and with a slight abuse of notation, write \(g\) as the realization
taking values in \(\Rbb^{d}\). Throughout, \(\sbound_g\) denotes the upper bound of the covariance spectral norm of \emph{any} random vector \(g\).
\begin{table}[h]
\centering
\caption{}
{\tiny
\begin{tabular}{|l|l|}
\hline
\((\tsigalg_t)_{t \in \Nz}\): Filtration of \(\tsigalg\) & \(\T_t\): test signal \\
\(\trace(M)\): trace of a matrix \(M\) & \(\ol{\athr}\): threshold (Algorithm \ref{alg:FTAAD}) \\
\(\opnorm{M}\): operator norm & \(\ol{d}_t\): martingale term (Lemma \ref{lem:a_sub_Gauss})\\
$\fnorm{M}$: Frobenius norm & \(q_t\):  predictable term (Lemma \ref{lem:predict_upperBnd}) \\
$\eig(M)$: eigenvalue of matrix \(M\) & \(\adver_t\): adversarial input \\
\(S^{\complift}\): set complement of \(S\) & \(M\) and \(\ol{M}\): quantities related to \(\athr\) (see \eqref{eq:aux_v1}) \\
\(\loro{t_1}{t_2}\): open interval & \(\sbound_0\): \(\sbound_g\) of the initial state \(\st_0\)\\
\hline
\end{tabular}
}
\end{table}
\section{Preliminaries}
\label{sec:ls}

\subsection{System description and problem statement}
\label{subsec:ps}
Consider a cyber-physical system (CPS) modeled as a linear
stochastic discrete-time process  
\begin{align}
\label{eq:lsys1}
        \stt{}{t+1} = A \stt{}{t} + B \cont{}{t} + \pnoise_{t+1} \text{ for \(t \in \Nz\),}
\end{align}
where \(\stt{}{t}\) and \(\cont{}{t}\) are the vectors of states and the control actions at time \(t\), taking values in the sets \(\stsp \subset \Rbb^d\) and \(\admact \subset \Rbb^m\), respectively. The matrices \(A \in \Rbb^{d \times d}\) and \(B \in \Rbb^{d \times m}\) denote the system and the control matrix, respectively. The sequence \(\pnoise \Let (\pnoise_t)_{t \in \Nz}\) is the process noise, which takes values in the set \(\nsp \subset \Rbb^d\) for each \(t\), and \(\st_0\) is the initial state of the system with a known distribution \(\mu_0\), and is assumed to be independent of the process noise  \(\pnoise\).
\begin{assum}
\label{assum:sub_gaussian_noise}
We enforce the following assumptions:
\begin{enumerate}[label=\textup{(\ref{assum:sub_gaussian_noise}-\alph*)}, leftmargin=*, widest=b, align=left]
    \item \label{it:sgn}\((\pnoise_t)_{t \in \Nz} \subset L^2(\ps)\) is an i.i.d sequence of \(\Rbb^d\mbox{-}\)valued sub-Gaussian random vectors with zero mean and covariance \(\var_{\pnoise} = \var_{\pnoise}^{\top}\succ 0\). 
   
    \item In addition to the data associated to \eqref{eq:lsys1}, \(\st_0\) is a zero mean sub-Gaussian random vector taking values in \(\stsp\) with covariance \(\Var(\st_0)\).
\end{enumerate}
\end{assum}
 Let  \(\policy \Let (\policy_t)_{t \in \Nz}\) be a deterministic Markov policy, where \(\policy_t: \Rbb^{d} \lra \Rbb^m\) for each \(t\), is a Borel measurable function.  We assume that the CPS \eqref{eq:lsys1} plays linear feedback control actions (generated by \(\policy_t\)) at each time \(t\), which are expressible as
\begin{align*}
    \cont{}{t} = \policy_t(\st_t) = \fb\st_t,
\end{align*}
for some state-feedback gain \(\fb \in \Rbb^{m \times d}\). The closed-loop dynamics is then given by
\begin{align}\label{eq:closed_lsys}
    \stt{}{t+1} &= (A  + B \fb) \st_t + \pnoise_{t+1}
    =  \Acl \st_t + \pnoise_{t+1},
\end{align}
where \(\Acl \in \Rbb^{d \times d}\) denotes the closed-loop system matrix. We further assume that the CPS operates under steady-state conditions with the process noise \(\pnoise\) acting as an external disturbance. Note that \((\cont{}{t})_{t \in \Nz}\) is a sequence of control actions that are optimally generated by solving an underlying optimal control problem with a specific cost function. The objective of the designed optimal control is to stabilize the system under a nominal condition. On the other hand, the adversarial attacker aims to tamper with the system's performance by either increasing the cost or causing the system to become unstable. We model the adversaries/attackers by an \(\Rbb^d\mbox{-}\)valued additive disturbance process \(\donoise{} \Let (\donoise{t})_{t \in \Nz} \subset \Rbb^d\) affecting the states in the following way:
\begin{equation}
    \label{eq:attack model}
    \sm{}{t} =  \stt{}{t} + \donoise{t} \quad \text{for }t \in \Nz,
\end{equation}
where \(\sm{}{t}\) is the vector of output measurements at time \(t\) and is available to the controller.

Under the preceding premise on the underlying CPS and the attackers, this article poses the question:
\begin{quote}
    \textsf{Q:} Given the stabilizable and the completely observable system \eqref{eq:lsys1} with the observation model \eqref{eq:attack model}, is it possible to detect, with high probability and in \emph{finite} time, \emph{additive adversarial attacks}? 
\end{quote}
We answer the posed question \textsf{Q} in the affirmative; see Algorithm \ref{alg:FTAAD} (AD-CPS) ahead in \S\ref{sec:mrg} and its subsequent analysis for details. In what follows, we will present AD-CPS and prove that the probability of incurring \emph{false positive} and \emph{false negative} errors by AD-CPS is small.

\subsection{Data available to defenders and attackers}\label{subsec:attack model}
\textit{Attackers: }
The adversaries may replace the true, honest observations with corrupt or fake data, focusing on attacks in which they inject an additive disturbance, represented by the process \(\donoise{}{}\) into the sensor channel.
%
The adversaries have access to the system parameters \((A,B)\), the observations \((\sm{}{t})_{t \in \Nz}\), and the detection algorithm adopted by designers. The actuators are assumed to be honest and the control policy \(\policy\) and the corresponding policy map (which is the state-feedback gain \(\fb\)) are inaccessible to adversaries. 
%
The intelligence of the attackers lies in how these fake or corrupt data are manufactured; in other words, the attacker may choose to inject an arbitrary process, or it may opt to lie dormant and undetected, subsequently injecting a crafted history-dependent adversarial disturbance to gradually deteriorate the performance of \eqref{eq:lsys1}.

\textit{Defenders: }The defenders have no access to the attack policies or to the mechanism by which they are generated. The only information available to them is the data produced by the compromised system. Moreover, the distribution of the process noise is unknown apart from an upper bound on the spectral radius of the process noise covariance.
\section{Attack detection algorithm for cyber-physical Systems (AD-CPS): intuition and discussions}
\label{sec:mrg}
Here we establish AD-CPS for \emph{additive adversarial attacks} as defined in \eqref{eq:attack model}. Recall from \S\ref{sec:intro} and \eqref{eq:attack model} that
\(\sm{}{t}\) observed at time \(t\) is termed as \emph{honest} if the \emph{measurement is honest}, and corrupt if the \emph{measurement is corrupt}.
In this setting, the question \textsf{Q} posed in \S\ref{sec:ls} falls within the broader ambit of \emph{online anomaly detection} under an adversarial setting. 
%
\subsection{Some preliminary results}\label{subsec:preliminary}
%
\begin{defn}
    \label{def:a_p_process}
     Let \(\bigl(\tsigalg_t \bigr)_{t \in \Nz}\) be a filtration of \(\tsigalg\). We say a stochastic process \(g = (g_t)_{t \in \Nz}\) is \emph{adapted} to \(\bigl(\tsigalg_t \bigr)_{t \in \Nz}\) if for each \(t \in \Nz\), \(g_t\) is \(\tsigalg_t\mbox{-}\)measurable. Moreover, \(g=(g_t)_{t \in \Nz}\) is said to be \(\bigl(\tsigalg_t \bigr)_{t \in \Nz}\) \emph{predictable} if for each \(t \in \Nz\), \(g_t\) is \(\tsigalg_{t-1}\mbox{-}\)measurable (with \(\tsigalg_{-1} \Let \aset[\big]{\emptyset,\ps}\)).
    \end{defn} 
The following result is central:
\begin{thm}[Doob decomposition \cite{ref:DW-91}]\label{th:doob_decom}
Let \(g =(g_t)_{t \in \Nz}\) be an \(\Rbb^d\mbox{-}\)valued random process adapted to \((\tsigalg_t)_{t \in \Nz}\), and \(g_t \in L^1(\ps)\) for each \(t \in \Nz\). Then the process \(g\) admits a Doob decomposition given by 
\begin{align}
    \label{eq:doob_decom}
    g_t = \mart_t + \predic_t + g_{t'}  \text{ for }t \in \Nz, \, \PP\mbox{-}\text{almost surely, where}
\end{align}
\begin{enumerate}[label=\textup{(\ref{eq:doob_decom}-\alph*)}, leftmargin=*, widest=b, align=left]
    \item \label{it:doob_m}\( (\mart_t)_{t \in \Nz}\) is a zero mean \((\tsigalg_t)_{t \in \Nz}\mbox{-}\)martingale with \(\mart_{t'} = 0\), and \(\mart_t\) at time \(t\) given by 
    \(
        \mart_t = \sum_{i=t'+ 1}^t\Bigl(g_i - \EE \cexpecof[\big]{g_i\given g_0,g_1,\ldots,g_{i-1}} \Bigr);
    \)

    \item \label{it:doob_p} \((\predic_t)_{t \in \Nz}\) is a \((\tsigalg_t)_{t \in \Nz}\mbox{-}\)predictable process  with \(\predic_{t'} = 0\) and 
    \(
        \predic_t = \sum_{i=t'+ 1}^t \Bigl(\EE \cexpecof[\big]{g_i\given g_0,g_1,\ldots,g_{i-1}} \Bigr) - g_{i-1} \Bigr).
    \)
\end{enumerate}
\end{thm}
\begin{rem}
    Note that in Theorem \ref{th:doob_decom} if \(t' = 0\) we end up with the conventional Doob decomposition \cite[\S \(12.11\)]{ref:DW-91}. However, over a fix window \(\window\), if \(t' = t-\window \), then \eqref{eq:doob_decom} is expressed as
    \(
        g_t = \sum_{i=t - \window +1}^t  \Bigl(g_i - \EE \cexpecof[\big]{g_i\given g_0,g_1,\ldots,g_{i-1}} \Bigr) + \sum_{i=t - \window+1}^t \Bigl(\EE \cexpecof[\big]{g_i\given g_0,g_1,\ldots,g_{i-1}} \Bigr) - g_{i-1} \Bigr) + g_{t-\window } \)  \(\PP\)-almost surely, for each \(t \in \Nz\).
\end{rem}


\subsection{Key intuition behind the detection test}\label{subsec:intuition}
This article establishes a class of anomaly detection-based algorithms where certain \emph{test signals} denoted by \(\T_t\), are tested against a suitably chosen threshold that depends on the upper-bound of the spectral norm of process noise covariance. More specifically, the primary idea of these classes of tests emerges from set memberships of the form 
\begin{equation}\label{eq:governing idea}
    \T_t \in \Gset_t \; \text{holds with \emph{high probability} for }t \in \Nz,
\end{equation} 
when the CPS is operating under nominal conditions. Here \(\Gset_t\) is appropriately chosen to minimize the error incurred by our algorithm (AD-CPS), and it depends on the statistical properties of the process noise. 
It will be revealed later in this section, we consider the case when the set membership is induced by a \emph{norm-based} test:
\begin{equation}\label{eq:blurb}
    \Gset_t \Let \aset[\big]{\T_t \suchthat \norm{\T_t} \leq \athr}.
\end{equation} 

Let us recall the control system \eqref{eq:closed_lsys} and suppose that \(\st_0, \st_1,\st_2,\ldots\) be a sequence of \(\Rbb^d\mbox{-}\)valued random vectors defined on \(\bigl(\ps,\Borelsigalg(\ps),\PP\bigr)\). Define  \(\bigl(\tsigalg_t \bigr)_{t \in \Nz}\) with \(\tsigalg_t \Let \sigma\bigl(\st_0,\ldots, \st_t\bigr)\) and \(\tsigalg_t \subset \tsigalg_{t+1}\) for every \(t \in \Nz\), and the filtered probability space \(\Bigl(\ps,\tsigalg,\bigl(\tsigalg_t \bigr)_{t \in \Nz},\PP \Bigr)\). Let \(\phist_t \Let \bigl(\stt{}{0},\ldots,\stt{}{t} \bigr)\) be the history until time \(t\).

Observe from Definition \ref{def:a_p_process} that \((\stt{}{t})_{t \in \Nz}\) is adapted to \((\tsigalg_t)_{t \in \Nz}\) (i.e., for each \(t \in \Nz\), \(\st_t\) is \(\tsigalg_t\mbox{-}\)measurable), and under Assumption \ref{it:sgn}, \(\EE \expecof[\big]{\norm{\stt{}{t}}}\) exists for all \(t \in \Nz\). The conditional mean of \(\stt{}{t+1}\) at time \(t\) given \(\phist_t\) is computed as
\(
    \EE \cexpecof[\big]{\stt{}{t+1}\given \phist_t} =  \Acl \EE \cexpecof[\big]{\stt{}{t}\given \phist_t} +  \EE \cexpecof[\big]{\pnoise_{t+1}\given \phist_t} = \Acl \stt{}{t}, 
\)
which implies that \(\stt{}{t+1} - \EE \cexpecof[\big]{\stt{}{t+1}\given \phist_t} = \pnoise_{t+1}\) for each \(t\). Similarly, we have \(\EE \cexpecof[\big]{\stt{}{t+1}\given \phist_t} - \stt{}{t} = \Acl \stt{}{t} - \stt{}{t}.\) 
Invoking Theorem \ref{th:doob_decom} and from \ref{it:doob_m}-\ref{it:doob_p}, \(\stt{}{t}\) is expressed as 
\begin{align}\label{eq:decomposition}
    \stt{}{t} &= \sum_{i=t-1}^{t} \bigl(\stt{}{i} - \EE \cexpecof[\big]{\stt{}{i}\given \phist_{i-1}}\bigr) + \sum_{i=t-1}^t \bigl(\EE \cexpecof[\big]{\stt{}{i}\given \phist_{i-1}} - \stt{}{i-1}\bigr) + \stt{}{t-2}\nn\\
    & \stackrel{\mathclap{(\dagger)}}{=} \sum_{i=t-1}^t \pnoise_{i} + \sum_{i=t-1}^t \bigl( \Acl \stt{}{i-1} - \stt{}{i-1
    }\bigr) + \stt{}{t-2} \quad \PP\mbox{-}\text{almost surely},\nn \\
    & \teL \mart_t + \predic_t + \st_{t-2} \quad \PP\mbox{-}\text{almost surely},
\end{align}
where the first term in \((\dagger)\) is the \emph{zero mean martingale part} (denoted by the process \(\mart\)) and the second term in \((\dagger)\) is the \emph{predictable part} (denoted by the process \(\predic\)) which depends on the system parameters \(\Acl\) and the past samples \((\st_{t-2},\st_{t-1})\). The decomposition in \eqref{eq:doob_decom} is carried out to extract the martingale \(\mart = (\mart_t)_{t \in \Nz}\), which corresponds to the cumulative process noise from the state process \((\st_t)_{t \in \Nz}\). Given the knowledge of certain statistical properties of the process noise, one can then define the sequence \((\Gset_t)_{t \in \Nz}\) and the threshold appropriately to execute the detection test. 
\begin{rem}
Note that, 
\(
    \st_{t-2} = \Acl^{t-2 } \st_0 + \sum_{j=0}^{t-3} \Acl^{t-3-j} \pnoise_{j+1}, 
\)
and it is a zero mean \(d\text{-}\)dimensional sub-Gaussian random vector with \(\norm{\var(\st_{t-2 })} \leq \sbound_{\st_{t-2}}\), where \(\sbound_{\st_{t-2 }} = \sbound_0 \norm{\Acl^{t-2 }}^2 + \vbound \sum_{j=0}^{t-3} \norm{\Acl^{t-3-j}}^2\).
Calculating \(\sbound_{\st_t}\) may be difficult, in general. One may compute \(\sbound_{\st_t}\) for any arbitrary \(t  \geq 1\) by the following iterations: for \(\ol{A}_0 = \identity{d}\),
\begin{align}
    \label{eq:iterations}
    \begin{cases}
        \ol{A}_t = \ol{A}_{t-1} \Acl,\, c_t = c_{t-1} + \norm{\ol{A}_{t-1}}^2   \text{ where }c_0 = 0,\\
        \sbound_{\st_t} = \sbound_0 \norm{\ol{A}_t}^2 + \vbound c_t \text{ where }\sbound_{\st_0} = \sbound_0.
    \end{cases}
\end{align}
\end{rem}
\subsection{The algorithm}\label{subsec:alg steps}
Recall the closed-loop CPS \eqref{eq:closed_lsys}. Within the context of our setting, we \emph{pick} the test signal \(\T_t\) in \eqref{eq:blurb} as
\begin{align}
    \label{eq:real_}
    \T_{t} &\Let \frac{1}{2}\sum_{\textcolor{blue}{s}=t-1}^t \Acl^{t-s}\sm{}{s} - \sm{}{t}- \frac{1}{2}\sum_{s=t-1}^t (\Acl^{t-s} - \identity{d})\sm{}{t-2}.
\end{align}
Using \eqref{eq:attack model} and \eqref{eq:decomposition}, the test signal  \eqref{eq:real_} can be expanded beginning from \(t-2\), and \emph{under the attack condition}, it admits the final expression 
\begin{align}
\label{eq:real}
    \T_{t} &=  \frac{1}{2}\sum_{s =t-1}^t\Acl^{t-s}\mart_s - \mart_t + \frac{1}{2}\sum_{s =t-1}^t\Acl^{t-s}\predic_s - \predic_t  + \adver_t, 
\end{align}
where \(\adver_t = \frac{1}{2}\sum_{s =t-1}^t\Acl^{t-s}\donoise{s} - \donoise{t} - \frac{1}{2}\sum_{s=t-1}^t (\Acl^{t-s} - \identity{d})\donoise{t-2}\) denotes the cumulative adversarial input. Note that under nominal operation, \((\st_t)_{t \in \Nz}\) is observed, and consequently, the anomaly \(\adver_t = 0\) for each iteration \(t\). However, when the system is under attack, \((\sm{}{t})_{t \in \Nz}\) is observed, which is used to construct the test signal \((\T_t)_{t \in \Nz}\), as defined in \eqref{eq:real_}. Here the term \(\ol{d}_t \Let \frac{1}{2}\sum_{s =t-1}^t\Acl^{t-s}\mart_s - \mart_t\) denotes the cumulative martingale at time \(t\), and \(q_t \Let \frac{1}{2}\sum_{s =t-1}^t\Acl^{t-s}\predic_s - \predic_t\) is the cumulative predictable term at time \(t\).

We now define the threshold to be employed in AD-CPS. Define the quantities
\begin{align}
\label{eq:aux_v1}
    \begin{cases}
        \vbound \Let \norm{\var_{\pnoise}} \text{ and  }\sbound_0 \Let \norm{\Var(\st_0)}, \\
        M \Let \Bigl(2 + 2 \norm{\Acl} + \norm{\Acl}^2 \Bigr) \quad \sbound_{d_t} \Let \vbound  M,\\
        \ol{M} \Let \biggl(4 + \frac{M}{4} + \norm{\Acl} \biggr) \quad \ol{\sbound} \Let \vbound\ol{M}.
    \end{cases}
\end{align} 
\begin{defn}
    \label{def:threshold}
    Fix a parameter \(\delta \in \loro{0}{1}\), we let \(\athr\) in Line \(3\) of Algorithm \ref{alg:FTAAD} be
    \begin{align*} 
        \athr \Let  \Bigl(\sqrt{2} + \sqrt{\ol{M}}\Bigr) \sqrt{k \vbound d \ln(1 / \delta)}, 
    \end{align*}
    where \(k\) is a constant dependent on the statistical properties of the process noise; we provide an estimate of \(k\) in the proof of Lemma \ref{lem:concen_ineq} in Appendix \ref{appen:proofs}.
\end{defn}
As mentioned in \S~\ref{subsec:intuition}, we consider threshold-based test by picking the set \(\Gset_t\) in the following way:
\[\Gset_t = \aset[\Big]{z \suchthat \norm{z} \leq \Bigl(\sqrt{2} + \sqrt{\ol{M}}\Bigr) \sqrt{k \vbound d \ln(1 / \delta)}}\; \text{for each }t.\]
\emph{We draw attention to the fact that the choice of \(\T \Let (\T_t)_{t \in \Nz}\); it is selected so that certain concentration-based results can be leveraged to establish strong theoretical guarantees for the errors. However, the user may suitably choose their own test signal \(\T\) to design custom tests.} In short, we provide a simple and effective framework, governed by \eqref{eq:governing idea}, for inferring the presence of adversaries.
\begin{algorithm}[!htb]
\caption{AD-CPS: An attack detection algorithm for cyber-physical systems}\label{alg:FTAAD}
\begin{algorithmic}[1]
\State {\bf{Data}}: flag = \texttt{null}, tuning parameter \(\tgain\);
\State {\bf{Initialize}}: Assume that \(\donoise{1} = 0\), to initialize the algorithm, calibrate \(\tgain\) and choose \(\ol{\athr} = \tgain \athr\);
\While{\( t \geqslant 2\)}
\State Observe \(\sm{}{t} \in \Rbb^{d}\)
\If{\(\norm{\T_t} \leq \ol{\athr}\)} \Comment{\(\T_t\) is defined in \eqref{eq:real}}
\State \texttt{Declare} \(\sm{}{t} \rightarrow\) \texttt{honest} 
\State \texttt{return} flag = 0 \Comment{Conclude that \(\sm{}{t}\) is honest}
\ElsIf{\(\norm{\T_t} > \ol{\athr}\)}
\State \texttt{Declare} \(\sm{}{t} \rightarrow\) \texttt{corrupt}  
\State \texttt{return} flag = 1  \Comment{Conclude that \(\sm{}{t}\) is corrupt}
\EndIf
\State \(t \gets t+1\)
\EndWhile
\end{algorithmic}
\end{algorithm}

We amplify several technical aspects of Algorithm \ref{alg:FTAAD} with the help of the following remarks:
%
%
\begin{rem}\label{rem:g_rem}
    Algorithm \ref{alg:FTAAD} provides a verifiable test to indicate, with high probability, whether \(\sm{}{t}\) observed at time \(t\), is honest or corrupt. The threshold \(\athr\) is chosen appropriately so that AD-CPS incurs fewer errors --- both false positive and false negative errors, theoretically --- with high probability (under certain mild assumptions); see \S\ref{sec:erros} and \S\ref{sec:fn_err} for more detail. However, we emphasize two points: first, minimizing both errors simultaneously is generally difficult, and consequently, a practitioner must make a trade-off when selecting \(\athr\);and second, in most practical situations, adjusting this trade-off requires appropriately tuning the theoretically obtained threshold in Definition \ref{def:threshold}. We acknowledge that this tuning is not optimal; see \S\ref{sec:ne} for further information on how the tuning is performed.
\end{rem}
\begin{rem}
    \label{rem:calibration}
    Offline tuning of the threshold is carried out by varying the tuning parameter \(\tgain\) (see Algorithm \ref{alg:FTAAD}). It is adjusted by computing the false positive error under nominal condition (i.e., when \(\donoise{t}=0\) for each \(t\)) as a function of \(\vbound\), for different values of \(\tgain\); see \S\ref{sec:ne} for more details.
\end{rem}
%

\section{Guarantees on false positive error}
\label{sec:erros}
In this section we establish that under the non-attack condition (i.e., when \(\donoise{t}=0\)) the \emph{false positive error} (FPE) incurred by AD-CPS is small for each \(t \in \Nz\) such that \(\sm{}{t}\) satisfies the conditions given in Definition \ref{def:good sample} (ahead). To that end, we define the two hypotheses:
\begin{align}
    \label{eq:hypotheses}
    \begin{aligned}
        \hypothesis_{0,t}:\donoise{t} = 0 \text{ and }
        \hypothesis_{1,t}: \donoise{t} \neq 0.
    \end{aligned}
\end{align}
 \begin{defn}\label{def:good sample}
 The observed state \(\sm{}{t}\) at time \(t \in \Nz\) is said to be \emph{\(2\)\text{-}step honest} if \(\sm{}{t-2}\) and  \(\sm{}{t-1}\) are honest, i.e., \(\donoise{t-2}= \donoise{t-1}=0\). 
 \end{defn}
The preceding definition is crucial for Theorem \ref{th:error_incurred}. The guarantee on FPE provided in Theorem \ref{th:error_incurred} is meaningful only when \(\sm{}{t}\) is \(2\)\text{-}step honest (in the sense of Definition \ref{def:good sample}). This is because, the cumulative effect of the past adversarial disturbances \((\donoise{t-2}, \donoise{t-1})\) may cause the system to behave abnormally, raising the detection alarm even though the null hypothesis \(\hypothesis_{0,t}\) is true at time \(t\), and significantly increasing the incurred FPE. 
Moreover,  \(2\)-step honesty is not crucial as far as the algorithm's operation is concerned but it is required to establish the assertion of Theorem \ref{th:error_incurred}.

 We are now ready to present the first main result of this article. Suppose that \( \ol{h} \Let \norm{\Acl - \identity{d}}\Bigl(\sum_{s =t-1}^t \norm{\Acl^{t-s}}+1\Bigr).\)
 %
\begin{thm}
    \label{th:error_incurred}
    Consider the control system \eqref{eq:lsys1}/\eqref{eq:closed_lsys} along with its associated data, and let Assumption \ref{assum:sub_gaussian_noise} hold. Consider Algorithm \ref{alg:FTAAD}, and fix \(\delta \in \loro{0}{1}\) as in Definition \ref{def:threshold}. Define the event
    \begin{equation}
        \label{eq:good_set}
        \goodset_t \Let  \aset[\Bigg]{\Bigg \lVert\frac{1}{2}\sum_{s =t-1}^t  \Acl^{t-s}\mart_s -\mart_t \Bigg \rVert \leqslant \sqrt{k \ol{\sbound} d \ln{\bigl(1/\delta \bigr)}}}, 
    \end{equation}
    where \(\ol{\sbound}\) is defined in \eqref{eq:aux_v1}. Suppose that there exist \(\beta>0\), \(\gamma \in \lcro{0}{1}\), and a compact set \(K\) such that 
    \begin{align}\label{eq:stability assump}
    \begin{cases}
        \EE \cexpecof[\big]{\norm{\st_1} \given \st_0} \leq \gamma \norm{\st_0}\quad \text{for all }\st_0 \notin K, \; \text{and}\\
        \beta \Let \sup_{\st_0 \in K} \EE \cexpecof[\big]{\norm{\st_1} \given \st_0} < +\infty.
    \end{cases}
    \end{align}
    Then on \(\goodset_t\) and on the set \(\aset[\big]{\sm{}{t} \text{ is \(2\)-step honest}}\) (in the sense of Definition \ref{def:good sample}), given \(\hypothesis_{0,t}\), at time \(t\) Algorithm \ref{alg:FTAAD} declares \(\sm{}{t}\) to be corrupt with probability at most 
    \begin{align*}
        \delta + \frac{2 \ol{h} \sum_{i=t-1}^t \Bigl( \gamma^{i-1} \norm{\st_0} + \beta \sum_{k=0}^{i-2} \gamma^{i-2-k} \Bigr)}{\sqrt{2k \vbound d \ln(1/\delta)}}.
    \end{align*}
\end{thm}
%
%
The proof of Theorem \ref{th:error_incurred} has been relegated to Appendix \ref{appen:t1_error}. In what follows, we include certain remarks to elaborate on Theorem \ref{th:error_incurred} from various perspectives.
%
\begin{rem}\label{rem:stability rem}
The following observations highlight several aspects of negative drift-type conditions \eqref{eq:stability assump} and the existence of an admissible control that adheres to \eqref{eq:stability assump} below:
\begin{itemize}[leftmargin=*, label = \(\circ\)]
    \item \label{it:general remark} The stability of the CPS manifests via \eqref{eq:stability assump}, and it adversely affect the false positive error probability. In other words, the system must be well-behaved under nominal operation to ensure low FPE. These conditions resemble negative drift conditions, \cite[Theorem \(4.3\)]{ref:SPM-RLT-92} and \cite[Proposition \(1\)]{ref:DC-JL-14}, and conveys that the magnitude of the control action must be large for states with large magnitude so that \(\EE \cexpecof[\big]{\norm{\st_1} \given \st_0} \leq \gamma \norm{\st_0}\) for \(\st_0 \notin K\) is satisfied. 

    \item \label{it:a suff condi} Recall from \eqref{eq:closed_lsys} that \(\Acl = A+BK\). A set of sufficient conditions for \eqref{eq:stability assump} are given by the following: 
    \begin{align*}
        \begin{cases}
        \Acl  \text{ is Lyapunov stable and } \eig_{\max}(\Acl^{\top}\Acl) <1,\, \text{with }\\
        \gamma = \frac{1}{2}\sqrt{\eig_{\min}(\Acl^{\top}\Acl)},\,
        K \Let \aset[\Big]{\dummyx \in \Rbb^d \suchthat \norm{\dummyx} \leq \frac{\sqrt{2\trace(\var_{\pnoise})}}{\gamma}}, \\
        \beta = \max_{\st \in K} \sqrt{\st^{\top}\Acl^{\top} \Acl \st} + \sqrt{\trace(\var_{\pnoise})}.
        \end{cases}
    \end{align*}

    \item \label{it:fpe is small} A careful scrutiny reveals that on picking \(\st_0 = 0\), one can simplify the expression \(2 \ol{h} \sum_{i=t-1}^t \Bigl( \gamma^{i-1} \norm{\st_0} + \beta \sum_{k=0}^{i-2} \gamma^{i-2-k} \Bigr)\). Indeed, observe that for \(i \in \aset[]{t-1,t}\), the finite sum \(\sum_{k=0}^{i-2} \gamma^{i-2-k} \leq \frac{1}{1-\gamma}\), implying that \(2 \ol{h} \sum_{i=t-1}^t \Bigl( \gamma^{i-1} \norm{\st_0} + \beta \sum_{k=0}^{i-2} \gamma^{i-2-k} \Bigr) \leq 2 \ol{h} \sum_{i=t-1}^t \biggl(  \frac{\beta}{1-\gamma}\biggr) 
        = 4 \ol{h} \biggl(  \frac{\beta}{1-\gamma}\biggr).\)
    We immediately observe that when \(\vbound\) is small (and consequently, \(\sqrt{\trace(\var_{\pnoise})}\) is also small), one can effectively design a state-feedback controller (without constraints on the magnitude of control actions) such that the closed-loop response of the CPS \eqref{eq:lsys1} remains bounded and close to zero under the applied state-feedback. This implies that \(\beta\) (defined above in Remark \ref{it:a suff condi}) is small and possibly close to zero, consequently leading to a low false positive error.
\end{itemize}
\end{rem}
\begin{rem}
    \label{rem:2-step honest}
    We emphasize that this article focuses on establishing a \emph{simple} detection scheme (AD-CPS) with \emph{non-asymptotic} guarantees, imposing relatively few assumptions on the adversarial disturbance \((\donoise{t})_{t \in \Nz}\) thus widening the ambit of our results. The \(2\)-step honesty assumption arises naturally given the fact that the designer possesses no information on how \((\donoise{t})_{t \in \Nz}\) is generated, i.e., the joint distribution of \((\donoise{t-2},\donoise{t-1},\donoise{t})\) is unknown. Note that it may be possible to provide strong non-asymptotic guarantees if the designer possesses data on the \emph{specific attack model} of adversaries. 
\end{rem}
    \begin{rem}
        \label{rem:choice of q}
        Please note that in general one may consider that the observed state \(\sm{}{t}\) at time \(t\) is \emph{\(\window\)-step honest} by looking back over a \(\window\)-length window from the current time instant \(t\) with \(\window >2\). The algorithm and the theoretical guarantees will follow through modulo certain straightforward changes. We emphasize that the studied case with \(\window=2\) is an example of a relatively myopic horizon and the system/test reacts \emph{quickly} to adversarial changes. For moderately large \(\window > 2\), it is natural to expect that the system will react, relative to the same test as outlined above, in a more graceful manner. For large values of \(\window\), one expects the test to behave sluggishly, leading to significant delays in detecting violent adversaries operating near the threshold boundary. Moreover, larger values of \(\window\) will \emph{significantly increase the threshold} due to the large variance of \(m = (m_t)_{t\in \Nz}\) computed over \(\window\), further weakening the theoretical results. Consequently, the algorithm will overlook certain immediately harmful adversaries, significantly damaging the CPS. 
        %
\end{rem}
\begin{rem}[Assumption on the adversaries]
    \label{rem:assump_adver}
    It is worth pointing out that no \emph{structural assumptions} (in the form of lower-bound or upper-bound on the magnitude of the adversary) have been imposed on the class of attacks. Moreover, we allow the unwanted disturbance process \(\donoise{}\) injected by the attackers to
     be a sequence of deterministic but unknown vectors, functions of the prior history of the CPS, or vectors randomized according to some unknown attack policy adopted by the adversary. AD-CPS and Theorem \ref{th:error_incurred} remain agnostic to these assumptions.
\end{rem}
\section{Guarantees on false negative error}\label{sec:fn_err}
Here we establish a guarantee on the \emph{false negative error} (FNE) incurred by AD-CPS. Ascertaining a strong guarantee on FNE without additional structural assumptions is challenging in contrast to Theorem \ref{th:error_incurred} where \emph{no} stipulations were imposed on the adversaries. Without such assumptions, nothing meaningful can be asserted, and hence, in the sequel (in Theorem \ref{th:t2_error_incurred}) we stipulate certain growth conditions on the cumulative adversarial inputs: Recall from \eqref{eq:real} that
\begin{equation}
    \label{eq:adv_inp}
    \adver_t =  \frac{1}{2}\sum_{s =t-1}^t\Acl^{t-s}\donoise{s} - \donoise{t} -\frac{1}{2}\sum_{s=t-1}^t (\Acl^{t-s} - \identity{d})\donoise{t-2}. 
\end{equation}
Given that \(\sm{}{t}\) is corrupt, on the event 
\(
    \aset[\big]{\norm{\adver_t} \geq \thres_{t}},
\) 
a non-asymptotic guarantee on the FNE can be provided. See Theorem \ref{th:t2_error_incurred} below for more details.

Recall the definition of the sequence \((\sbound_{\st_{t}})_{t \in \Nz}\) from \eqref{eq:iterations} and the quantity \(\ol{h}\) from the hypothesis of Theorem \ref{th:error_incurred}. Fix \(\delta \in \loro{0}{1}\) as in Definition \ref{def:threshold},
    \begin{align}\label{eq:adv_lb}
        &\lbt_t(\Acl,d,\vbound,\delta)\nn \\&\Let  \sqrt{d \vbound \ol{M}} + \ol{h}\sum_{i=t-1}^t \Bigl( \sqrt{d \sbound_{\st_{i-1}}} + \sqrt{k\sbound_{\st_{i-1}} \ln{(1/\delta)}} \Bigr).
    \end{align}
\begin{assum}
    \label{assum:ad_assum}
    In keeping with the notations established in \S\ref{sec:mrg} we suppose that\footnote{In the sequel, we suppress the dependency of \(\thres_t\) on \((\Acl,d,\vbound,\delta)\) and equivalently write \(\thres_t  = \thres_{t}(\Acl,d,\vbound,\delta).\)}
    \begin{align*}
        \thres_{t} = \lbt_t(\Acl,d,\vbound,\delta) + \athr \; \text{ with }\norm{\adver_t} \geq \thres_{t}.  
    \end{align*} 
\end{assum}
Note that in Assumption \ref{assum:ad_assum} and Theorem \ref{th:t2_error_incurred} ahead, we do not vary \(\delta\), essentially fixing the attack class in Assumption \ref{assum:ad_assum}. The associated false negative error rate is computed for that fixed attack class.

Under Assumption \ref{assum:ad_assum}, the following result yields a certificate on the false negative error incurred by AD-CPS:
\begin{thm}
    \label{th:t2_error_incurred}
    Consider the control system \eqref{eq:lsys1}/\eqref{eq:closed_lsys} along with its associated data, and let Assumption \ref{assum:sub_gaussian_noise} hold. Consider Algorithm \ref{alg:FTAAD}, and fix \(\delta \in \loro{0}{1}\) as in Definition \ref{def:threshold}.
    Define \(\eps \Let \frac{\ol{h}}{2}\sum_{i=t-1}^t \sqrt{k \sbound_{\st_{i-1}} \ln{(1/\delta)}}\),
    and recall from \eqref{eq:adv_inp}, the definition of \(\adver_t\).
    Under Assumption \ref{assum:ad_assum}, i.e., \(\thres_{t} = \lbt_t(\Acl,d,\vbound,\delta) + \athr\), on \(\aset[\big]{\norm{\adver_t} \geq \thres_t}\) and given \(\hypothesis_{1,t}\), at time \(t\) Algorithm \ref{alg:FTAAD} declares \(\sm{}{t}\) to be honest with probability at most \(\exp{\Bigl( - \frac{\eps^2}{k \vbound \ol{M} }\Bigr)} + 2 \delta\), where \(\ol{M}\) is defined in \eqref{eq:aux_v1}. 
\end{thm}
We place the proof of Theorem \ref{th:t2_error_incurred} in Appendix \ref{appen:t2_error}. Here we discuss a few technical aspects of Theorem \ref{th:t2_error_incurred}.
    \begin{rem}\label{rem:tradeoff_fpe fne}
        It is well known that changing the threshold effectively alters the balance between false-positive and false-negative errors. In our context, the threshold is chosen to effectively reduce the false-positive error. Consequently, the corresponding trade-off is reflected in Assumption \ref{assum:ad_assum} in Theorem \ref{th:t2_error_incurred}, which states that the cumulative effect of the adversary must exceed the designed threshold by \(\lbt_t(\Acl,d,\vbound,\delta)\). 
        We restrict the class of attacks to those that satisfy Assumption \ref{assum:ad_assum}, which is a natural artifact of the trade-off. 
        The case of stealthy attackers that respect the threshold yet cause serious damage to the CPS is not investigated in this article and will be treated elsewhere. 
    \end{rem}
\begin{rem}
    \label{rem:technical FNE}
    The following remarks highlight several aspects of Assumption \ref{assum:ad_assum} placed on the adversary.
    \begin{itemize}[leftmargin=*, label = \(\circ\)]
        \item \label{rem:on k_jt} Since the closed-loop system matrix is Schur stable, we have \(\max \limits_{1\leq i \leq d}\abs{\eig_i(\Acl)} \leq 1\). If \(\Acl\) is diagonalizable, then one can obtain an upper-bound on \(\sbound_{\st_t}\):
    \begin{align*}
        \limsup_{t \ra +\infty} \sbound_{\st_{t}} &\leq \sbound_0 \limsup_{t \ra +\infty} \norm{\Acl^{t}}^2 \\&+ \vbound \limsup_{t \ra +\infty} \sum_{j=0}^{t-1} \Bigl(\max \limits_{1\leq i \leq d}\abs{\eig_i(\Acl)}\Bigr)^{2(t-1-j)}\\
        &  \xrightarrow[t \to + \infty]{} \frac{\vbound}{1- \Bigl(\max \limits_{1\leq i \leq d}\abs{\eig_i(\Acl)}\Bigr)^2}.
    \end{align*} 
    On the other hand, if \(\Acl\) is not diagonalizable (but the Jordan normal form exists), obtaining a closed-form expression of the upper bound of \(\sbound_{\st_{t}}\) is difficult. However, one can show that \(\norm{\Acl^t} \lra 0\) as \(t \lra + \infty,\) and \(\limsup \limits_{t \to + \infty}\frac{\norm{\Acl^{t+1}}}{\norm{\Acl^t}} < 1\) due to \(\Acl\) being Schur stable. Hence, invoking \cite[Theorem 3.34]{WR:64}, the series \(\sum_{j=0}^{t-1} \norm{\Acl^{t-1-j}}^2\) converges to a finite quantity, implying \(\limsup_{t\ra +\infty}\sum_{j=0}^{t-1} \norm{\Acl^{t-1-j}}^2 < +\infty\).

    \item \label{rem:on adv assumpt} As pointed out above, the sequence \((\sbound_{\st_t})_{t \in\Nz}\) is a bounded sequence. This immediately implies that the \(\bigl(\lbt_t(\Acl,d,\vbound,\delta)\bigr)_{t \in \Nz}\) is also bounded. In fact one can obtain both lower-bound and upper-bound: Define \(\nbound \Let \sup_{t \in \Nz} \sbound_{\st_t}\). For an arbitrary \(t \in \Nz\), we observe that \(\sqrt{d \vbound \ol{M}}+ \ol{h}\sum_{i=t-1}^t \Bigl( \sqrt{d \sbound_0} + \sqrt{k\sbound_0 \ln{(1/\delta)}} \Bigr)\leq \lbt_t(\Acl,d,\vbound,\delta) \leq \sqrt{d \vbound \ol{M}} + 2 \ol{h} \sqrt{k \nbound d \ln{(1/\delta)}}.\)
    When \(\Acl\) is diagonalizable, then \(\nbound \geq \frac{\vbound}{1- \max \limits_{1\leq i \leq d}\abs{\eig_i(\Acl)}}\), and the upper-bound on \(\lbt_t(\Acl,d,\vbound,\delta)\) computed above is tighter. Hence, it suffices to impose that \(\norm{\adver_t} \geq \sup_{t \in \Nz} \thres_t \geq \limsup_{t \in \Nz} \lbt_t(\Acl,d,\vbound,\delta) + \athr = \sqrt{d \vbound \ol{M} } + 2 \ol{h} \sqrt{k \nbound d \ln{(1/\delta)}} + \Bigl(\sqrt{2} + \sqrt{\ol{M}}\Bigr) \sqrt{k \vbound d \ln(1 / \delta)}\).
    and the guarantee in Theorem \ref{th:t2_error_incurred} holds.
    \end{itemize}
\end{rem}
\begin{rem}
    \label{rem:fne is small}
    Observe that Theorem \ref{th:t2_error_incurred} does not require \(\sm{}{t}\) to be \emph{2\text{-}step honest} in the sense of Definition \ref{def:good sample}.
    Moreover, the probability estimated in Theorem \ref{th:t2_error_incurred} is small and saturates at \(2 \delta\) uniformly over time; see Figure \ref{fig:err_fne} ahead in the \S \ref{sec:ne}. 
\end{rem}
\section{Numerical experiment}\label{sec:ne}
The goal of this section is to empirically analyze various aspects of AD-CPS with respect to the `evaluation metrics' defined below. For the experiments, we extend the setting considered in this article to the partially observed case and design a \emph{residual-based} detector to identify a class of deception attacks. We extensively study the features and shortcomings of AD-CPS and comment on each assumption, especially the issue of violating \(\window\)-step honesty in Definition \ref{def:good sample}, as well as other conditions considered in the article. Our experiments showed that for benign process noise, violating \(\window\)-step honesty, the false positive error does not necessarily increase the false positive error (see Figure \ref{fig:violation}), and therefore, the assumption is not overly restrictive. We also compare AD-CPS with the optimal watermarking CUSUM-based detector developed in \cite{ref:AN-AT-AA-SD-23}, and our experiments validate that AD-CPS achieves comparable performance; see Table \ref{tab:comparison} and Figure \ref{fig:comparison} for further details.
\subsection{Partially observed setting}\label{subfig:partially observed}
A partially observed setting is considered for the experiment, which is given by
\begin{align}
    \label{eq:partially observed}
    &\st_{t+1} = A \st_t + B \cont{}{t} + \pnoise_t \text{ given initial condition }\st_0  \nn\\
    &y_t = C \st_t + \ol{w}_t, \text{ for }t \in \Nz,
\end{align}
where \(y_t \in \Rbb^{p}\) is the output at time \(t\), and \(C \in \Rbb^{p \times d}\) is the output matrix. The process noise \(\pnoise \Let (\pnoise_t)_{t \in \Nz} \subset \Rbb^{d}\) and the measurement noise \(\ol{\pnoise}\Let (\ol{{\pnoise}_t})_{t \in \Nz} \subset \Rbb^{p}\) are sequences of i.i.d. zero mean  Gaussian random vectors with covariances \(\Sigma_{\pnoise}\) and \(\Sigma_{\ol{w}}\), respectively. Of course \(\norm{\Sigma_{\pnoise}} \le \vbound\) and \(\norm{\Sigma_{\ol{\pnoise}}} \le \ol{\sbound}_{\pnoise}\). A standard Kalman filter satisfying the recursions
\begin{align*}
    &\wh{\st}_{t+1|t+1} = A \wh{\st}_{t|t} + B \cont{}{t} + L \res_t \text{ with }\wh{\st}_{0|0} = 0\\
    & \res_t = y_t - C \wh{\st}_{t|t-1},
\end{align*}
is designed to stabilize the system. The gain \(L = PC^{\top}(C P C^{\top} + R){\inverse}\), where \(P \succ 0\) is the solution of an algebraic Ricatti equation; we refer the readers to \cite[\S2]{ref:AN-AT-AA-SD-23} for further information about the Kalman filter. The sequence \((\res_t)_{t \in \Nz}\) denotes the innovation sequence (also known as the residual).
We considered a residual-based test signal: 
\begin{align}
    \label{eq:residual-based test signal}
    \T_t = \frac{1}{\window} \sum_{i=t-\window+1}^{t} \res_i - \res_t \text{ for every }t \in \Nz,
\end{align}
to detect adversaries. Here \(\window>0\) is a finite running window. It is well-known that for the system \eqref{eq:partially observed} operating under nominal condition and Gaussian assumptions, the residual follows \(\res_t \sim \mathcal{N}(0, \Sigma_{\res})\) for every \(t\), where \(\Sigma_{\res} = C P C^{\top} + \Sigma_{\ol{\pnoise}}\). Let us stipulate that \(\norm{\Sigma_{\res}} \le \sbound_{\res}\).

\textbf{System description. }
The inverted pendulum \cite{ref:ip} on a movable cart with a restricted cart length is considered for the experiments. The inverted pendulum is pivoted to the movable cart in an upright position with an unstable equilibrium. 
The system parameters are
\begin{align}\label{eq:example_IP}
    A = \begin{pmatrix}
        1 & 0.01 & 0.00011 & 0\\
        0 & 0.9982 & 0.0267 & 0.0001\\
        0 & 0 & 1.0016 & 0.01\\
        0 & -0.0045 & 0.3119 & 1.0016
    \end{pmatrix}\quad B = \begin{pmatrix}
        0.0001\\0.0182\\0.0002\\0.0454
    \end{pmatrix},
\end{align}
and we considered \(C = \identity{d}\). The control input \(\cont{}{t}\) at time \(t\) denotes the horizontal force acting on the moving cart. 

\textbf{Evaluation metrics.} 
For all experiments we fixed \(\delta = 0.01\). The time horizon for experiments was picked as \(T = 1000\). Let \(T_{1}\) and \(T_{2}\) denote the time instants at which the adversary begins and ceases to act, respectively, and define \(\Delta \Let T_2 - T_1>0\) to be the time interval when the adversary is active. The following metrics were considered to assess the performance of AD-CPS for all experiments:
\begin{itemize}[leftmargin=*]
    \item \emph{False positive error}, empirically computed as
    \begin{equation*}
        \wh{FPE}_{\window, \delta} \Let \frac{\text{number of threshold violations}}{T - \Delta},
    \end{equation*}
    under \emph{nominal operation}.
    \item \emph{Detection rate}, defined by
    \begin{equation*}
        \wh{DR}_{\window,\delta} \Let \frac{\text{number of detections under attack}}{\Delta},
    \end{equation*}
    The empirical \emph{false negative error} is defined to be \(\wh{FNE}_{\window,\delta} \Let 1 - \wh{DR}_{\window,\delta}\).
\end{itemize}

\textbf{Attack model.} For all experiments, we considered the deception attack \cite{ref:AN-AT-AA-SD-23} described by the process
\begin{equation}
    \label{eq:deception attack}
    \donoise{t} = A_a \donoise{t-1} + \wh{\donoise{}}_t
\end{equation}
where \(A_a \in \Rbb^{p \times p}\), and \(\wh{\donoise{}}_t \sim \mathcal{N}(0, \Sigma_a)\) for every \(t\), refers to the uncertainties associated with the adversary. During an attack, an adversary, with appropriately chosen \((A_a, \Sigma_a)\), may either \emph{replace the true output signal} \(y\) or may \emph{inject} \(v\) into \(y\) such that the corrupt data available to the system is given by
\begin{equation*}
    z_t = y_t + v_t \text{ for }t \in \Nz.
\end{equation*}
Naturally, under attack, the residual signal ceases to be a Gaussian random vector. Let us assume that \(\norm{\Sigma_a} \le \sbound_a\).

\textbf{Parameters. } We fixed \(\vbound\) at \(0.001\) across all experiments. The parameters \((\sbound_{\ol{\pnoise}}, \sbound_a, \tgain)\) are varied to test the performance of AD-CPS under different scenarios. Throughout, we fix the window size to be \(20\), except in Figure \ref{fig:fpe_nominal}. The system is initialized randomly (seed equal to \(2\)).

\subsection*{Results and discussions }
\begin{figure}[h]
\centering
\begin{subfigure}{0.45\linewidth}
    \includegraphics[width = 1\columnwidth, height = 2.2cm]{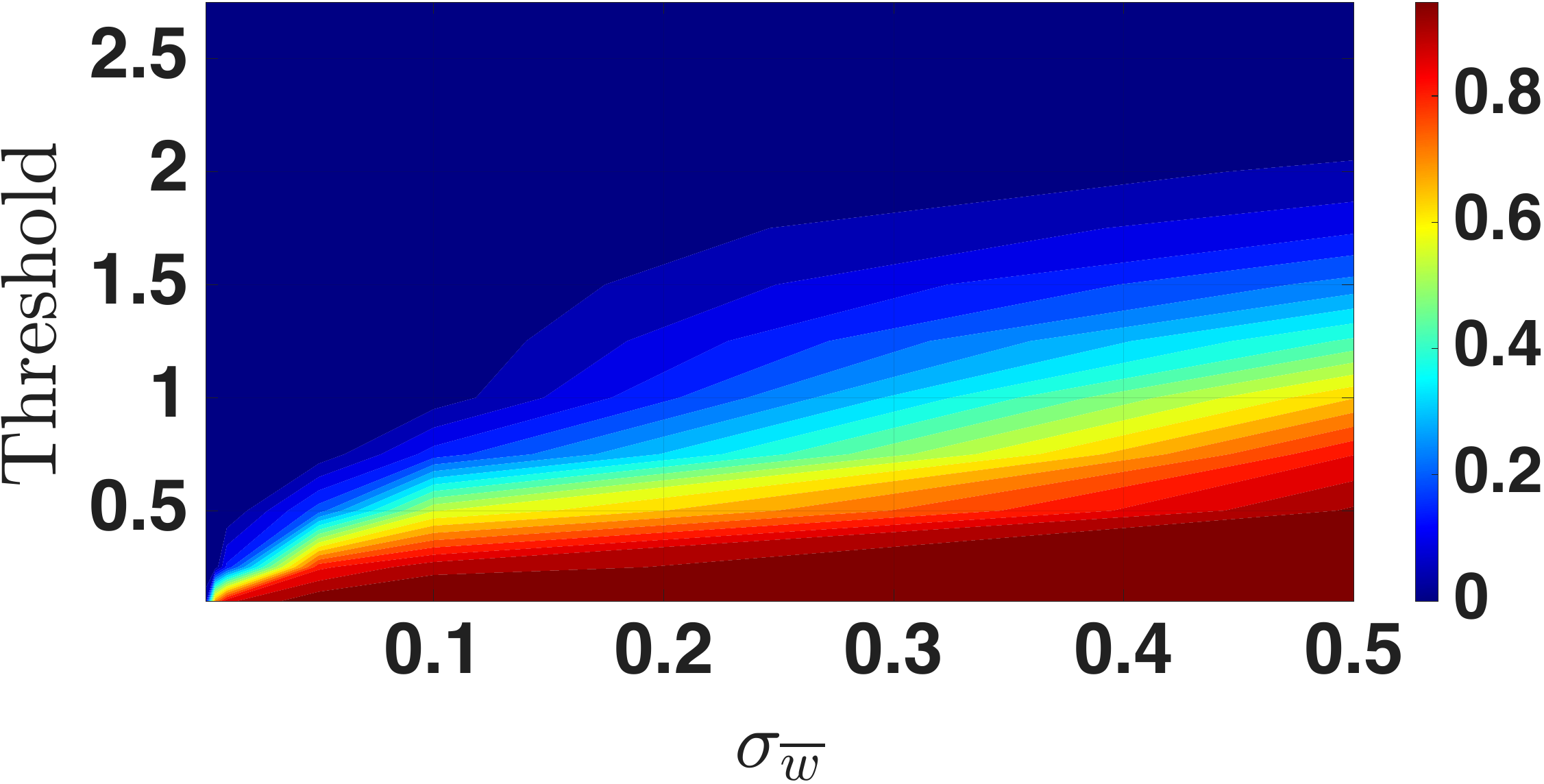}
    \caption{\(W = 5\)}
    \label{subfig:W_5}
\end{subfigure}
\begin{subfigure}{0.45\linewidth}
    \includegraphics[width = 1\columnwidth,height = 2.2cm]{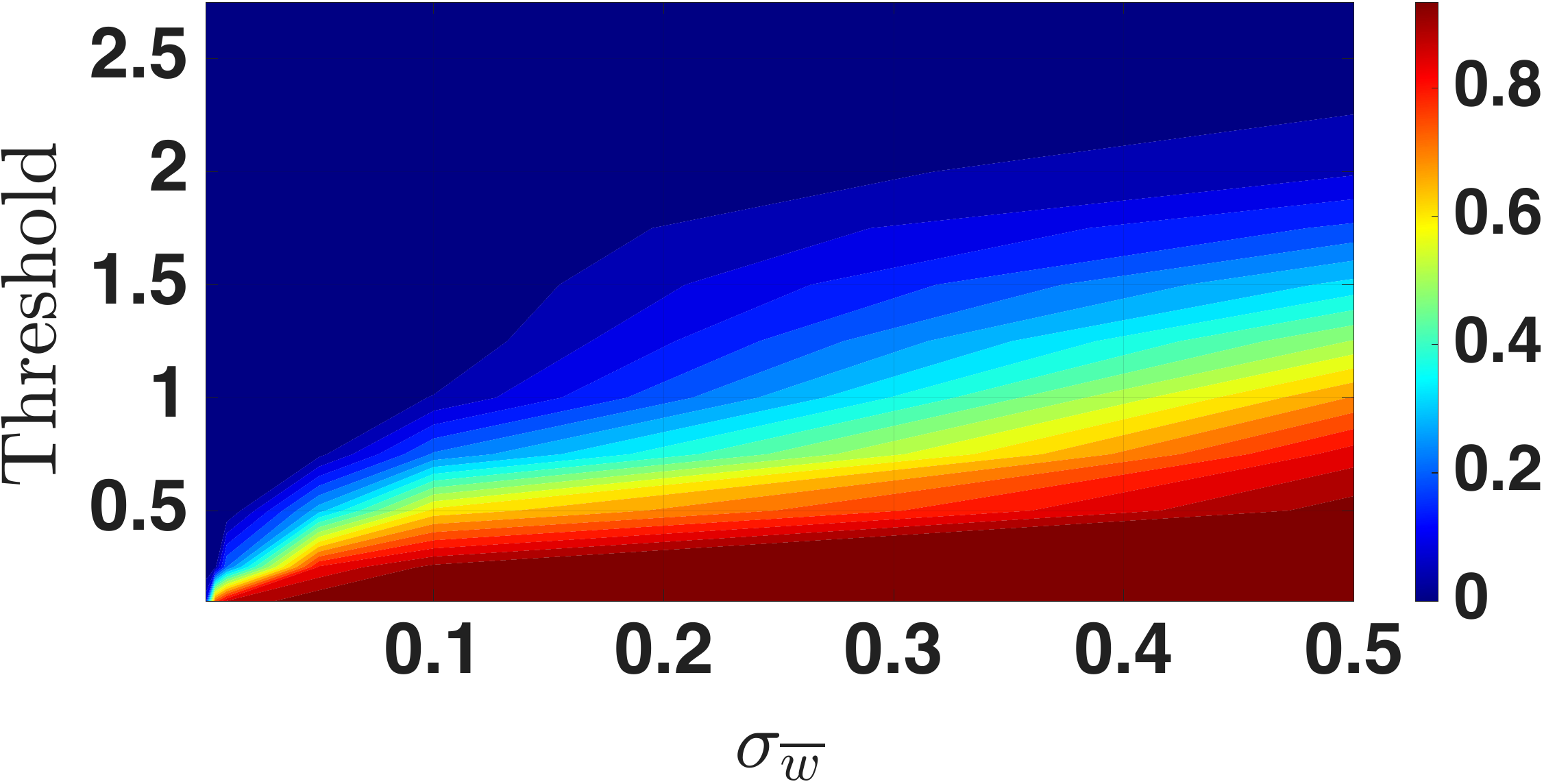}
    \caption{\(W=20\)}
    \label{subfig:W_20}
\end{subfigure}  
\vspace{-2mm}
\caption{Empirical false positive error \(\wh{FPE}_{\window,\delta}\) for \(\window = 5\) (Figure \ref{subfig:W_5}) and \(\window = 20\) (Figure \ref{subfig:W_20}) under nominal condition.}
\label{fig:fpe_nominal}
\end{figure}

\textbf{Tuning the threshold. }Figure \ref{fig:fpe_nominal} shows the variation of \(\wh{FPE}_{\window,\delta}\) as a function of the threshold \(\ol{\athr}\) and \(\sbound_{\ol{\pnoise}}\), under nominal condition. Two observations are noteworthy: first, for a fixed \(\sbound_{\ol{\pnoise}}\), the threshold must be increased to achieve a low false positive error. Second, for a fixed value of the threshold, as \(\sbound_{\ol{\pnoise}}\) is increased, \(\wh{FPE}_{\window,\delta}\) increases. Based on these observations, one may select, with reasonable accuracy, a suitable \(\tgain\) and threshold to attain a specific \(\wh{FPE}_{\window,\delta}\). 

\begin{figure}[ht]
\centering
\begin{subfigure}{0.45\linewidth}
    \includegraphics[width = 1\columnwidth, height = 2.2cm]{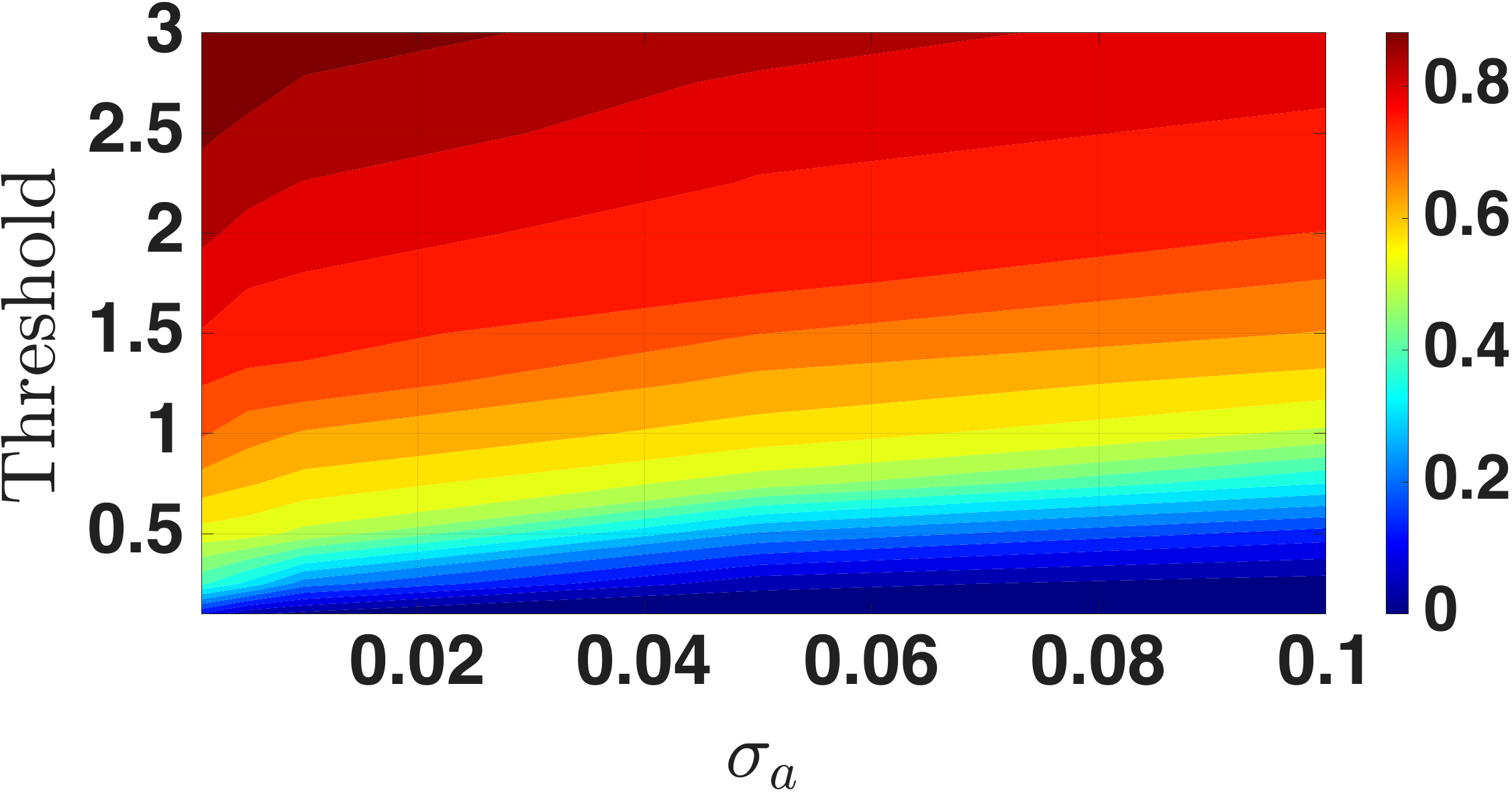}
    \caption{\(\sbound_{\ol{\pnoise}} = 0.005\)}
    \label{subfig:sigma_005}
\end{subfigure}
\begin{subfigure}{0.45\linewidth}
    \includegraphics[width = 1\columnwidth,height = 2.2cm]{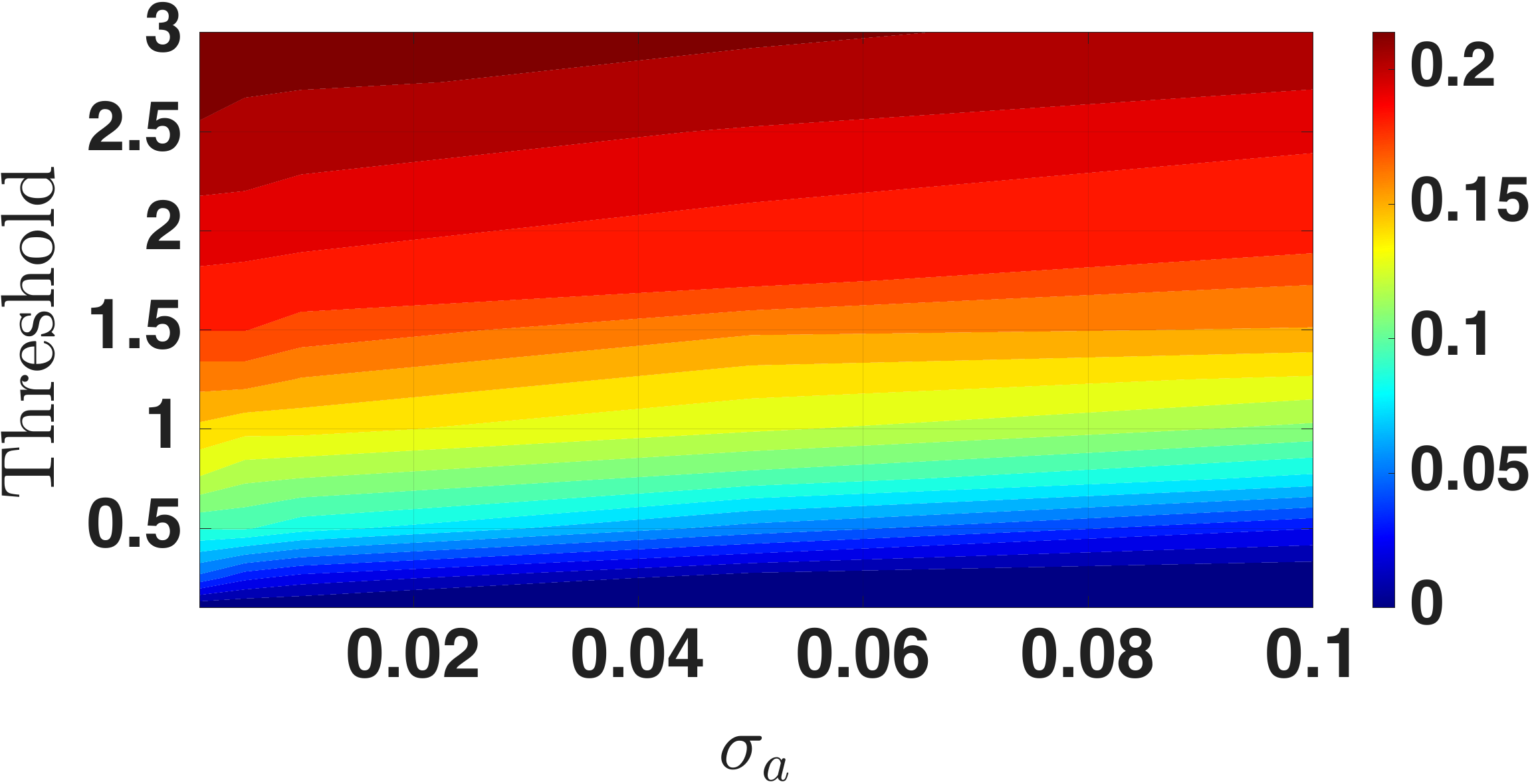}
    \caption{\(\sbound_{\ol{\pnoise}} = 0.1\)}
    \label{subfig:sigma_1}
\end{subfigure}  
\vspace{-2mm}
\caption{Empirical false negative error \(\wh{FNE}_{\window,\delta}\) for \(\sbound_{\ol{\pnoise}} = 0.005\) (Figure \ref{subfig:W_5}) and \(\sbound_{\ol{\pnoise}} = 0.1\) (Figure \ref{subfig:W_20}), as a function of \(\ol{\athr}\) and \(\sbound_{a}\). Here \(\window =20\) is considered.}
\label{fig:fne}
\end{figure}

\textbf{Performance. }Proper tuning is essential to improve the detection rate for systems under attack. It can be observed that for large \(\sbound_{\ol{\pnoise}}\), violating the conditions analogous to \eqref{eq:stability assump}, a high threshold must chosen to decrease \(\wh{FPE}_{\window,\delta}\). However, as demonstrated in Figure \ref{subfig:sigma_1}, if a high threshold \(\ol{\athr}\) is chosen, the \(\wh{FNE}_{\window,\delta}\) deteriorates drastically, even when the adversary employs benign attacks with small \(\sbound_a\). This observation immediately reveals that if the system is inherently unstable, our detector will not achieve good performance --- \(\wh{FNE}_{\window,\delta}\) increases as well --- highlighting that the conditions \eqref{eq:stability assump} are critical \emph{but not restrictive}. 

\begin{figure}[h]
    \centering
    \includegraphics[scale = 0.12]{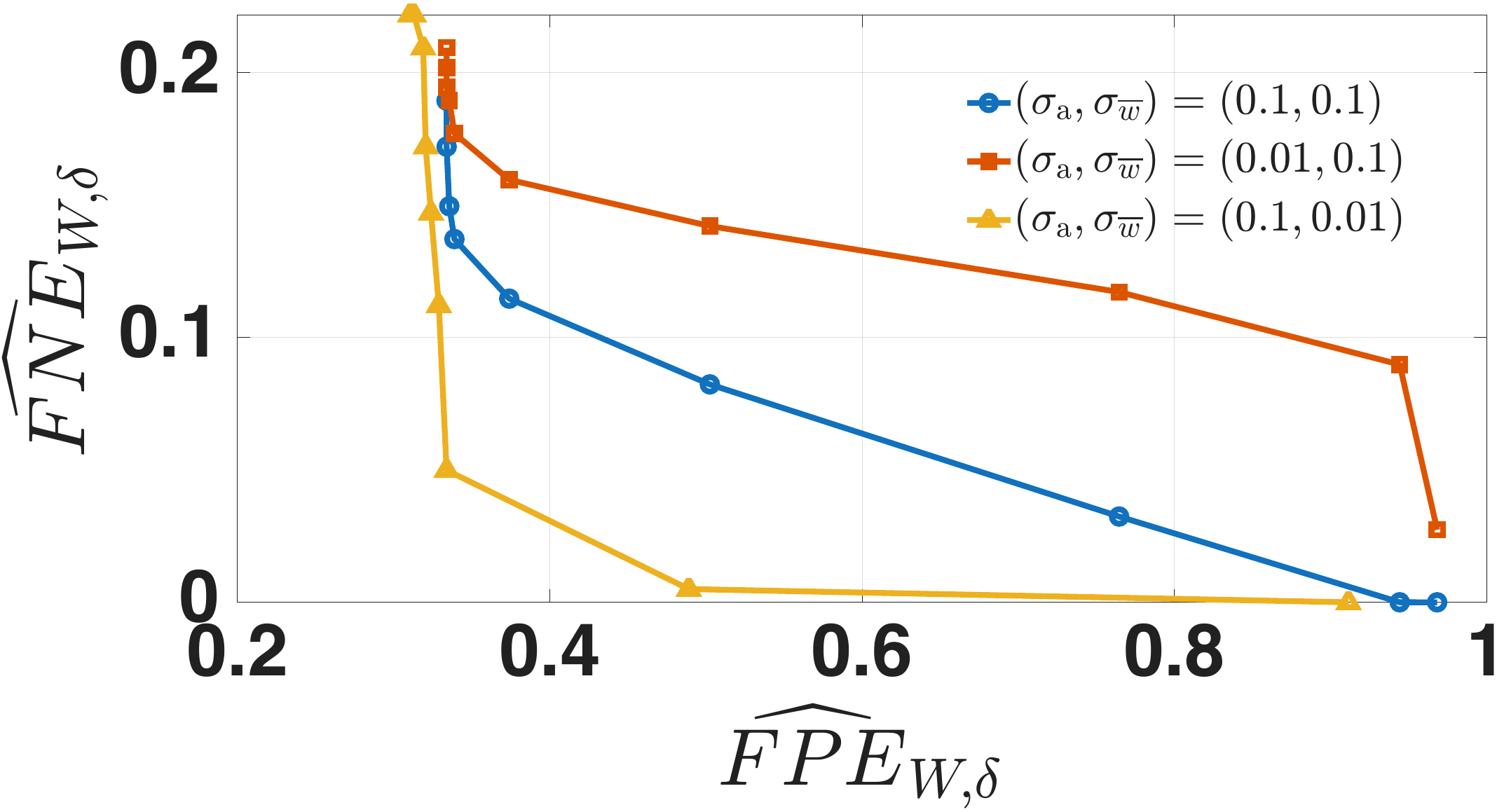}
    \caption{Trade-off between \(\wh{FPE}_{\window,\delta}\) and \(\wh{FNE}_{\window,\delta}\) for different values of \(\vbound\) and \(\sbound_{a}\).}
    \label{fig:trade_off}
\end{figure}

Figure \ref{fig:fne} summarizes the variation of \(\wh{FNE}_{\window,\delta}\) for a variable degree of stealthiness. If \(\sbound_{\ol{\pnoise}}\) is big, naturally, \(\wh{FNE}_{\window,\delta}\) would increase, as shown in Figure \ref{fig:fne}. The inherent trade-off between the false positive error and the false negative error, addressed in Remarks \ref{rem:g_rem} and \ref{rem:tradeoff_fpe fne}, is depicted in Figure \ref{fig:trade_off}. As threshold increase (from right to left), \(\wh{FPE}_{\window,\delta}\) decreases but \(\wh{FNE}_{\window,\delta}\) increases.

\begin{figure}
    \centering
    \includegraphics[width=0.55\linewidth, height = 2.4cm]{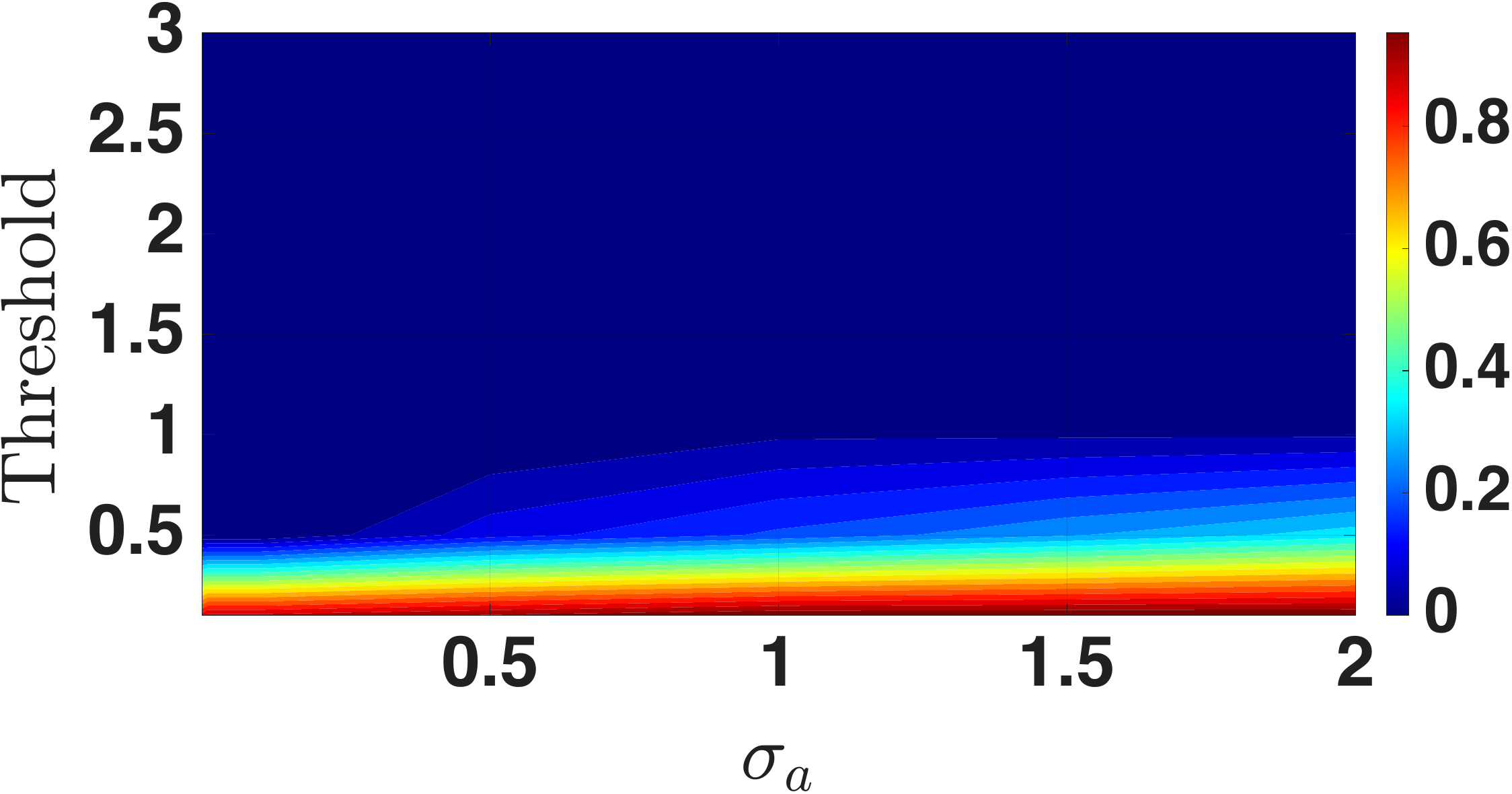}
    \caption{\(\wh{FPE}_{\window,\delta}\) when \(\window\)-step honesty is violated (Assumption \ref{assum:ad_assum}).}
    \label{fig:violation}
\end{figure}

Figure \ref{fig:violation} demonstrates the performance of AD-CPS when \(\window\)-step honesty, in Definition \ref{def:good sample}, is violated. When \(\sbound_a\) admits small values, \emph{even when \(W\)-step honesty is violated}, \(\wh{FPE}_{\window,\delta}\) does not deteriorates. This observation implies that the violation/falsity of \(W\)-step honesty \emph{does not} necessarily lead to high FPE, in general. Of course, the threshold must be chosen appropriately (from Figure \ref{fig:fpe_nominal}) to maintain a desired FPE tolerance. However, when \(\sbound_{\ol{\pnoise}}\) increases, \(\wh{FPE}_{\window,\delta}\) is relatively high. One can choose a higher threshold to effectively reduce \(\wh{FPE}_{\window,\delta}\), albeit at the expense of an increased false negative error, consistent with Figure \ref{fig:trade_off}.

\begin{figure}[ht]
\centering
\begin{subfigure}{\linewidth}
    \centering
    \begin{subfigure}{0.48\linewidth}
        \includegraphics[width=\linewidth,height = 2.7cm]{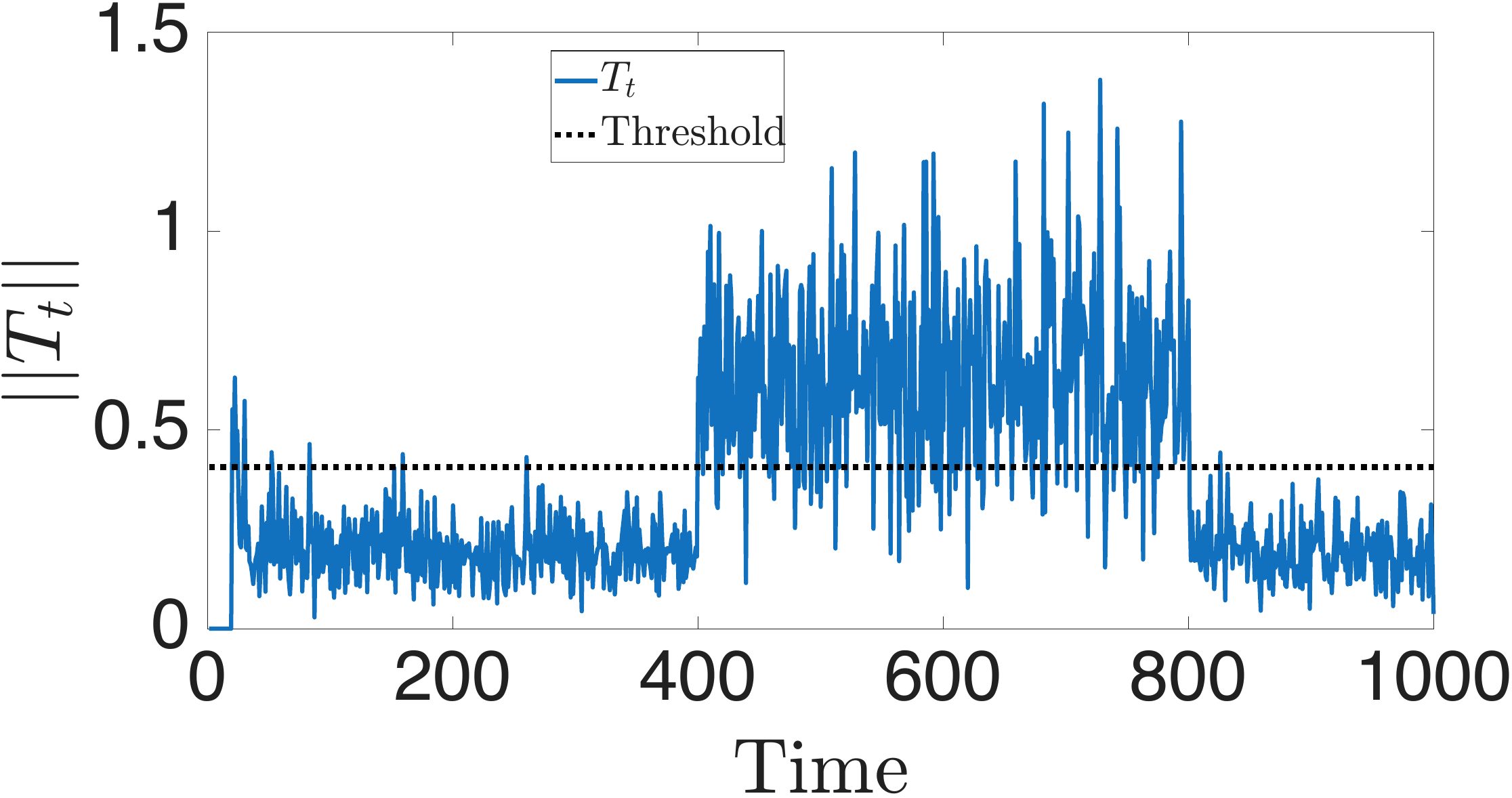}
    \end{subfigure}
    \hfill
    \begin{subfigure}{0.48\linewidth}
        \includegraphics[width=\linewidth,height = 2.56cm]{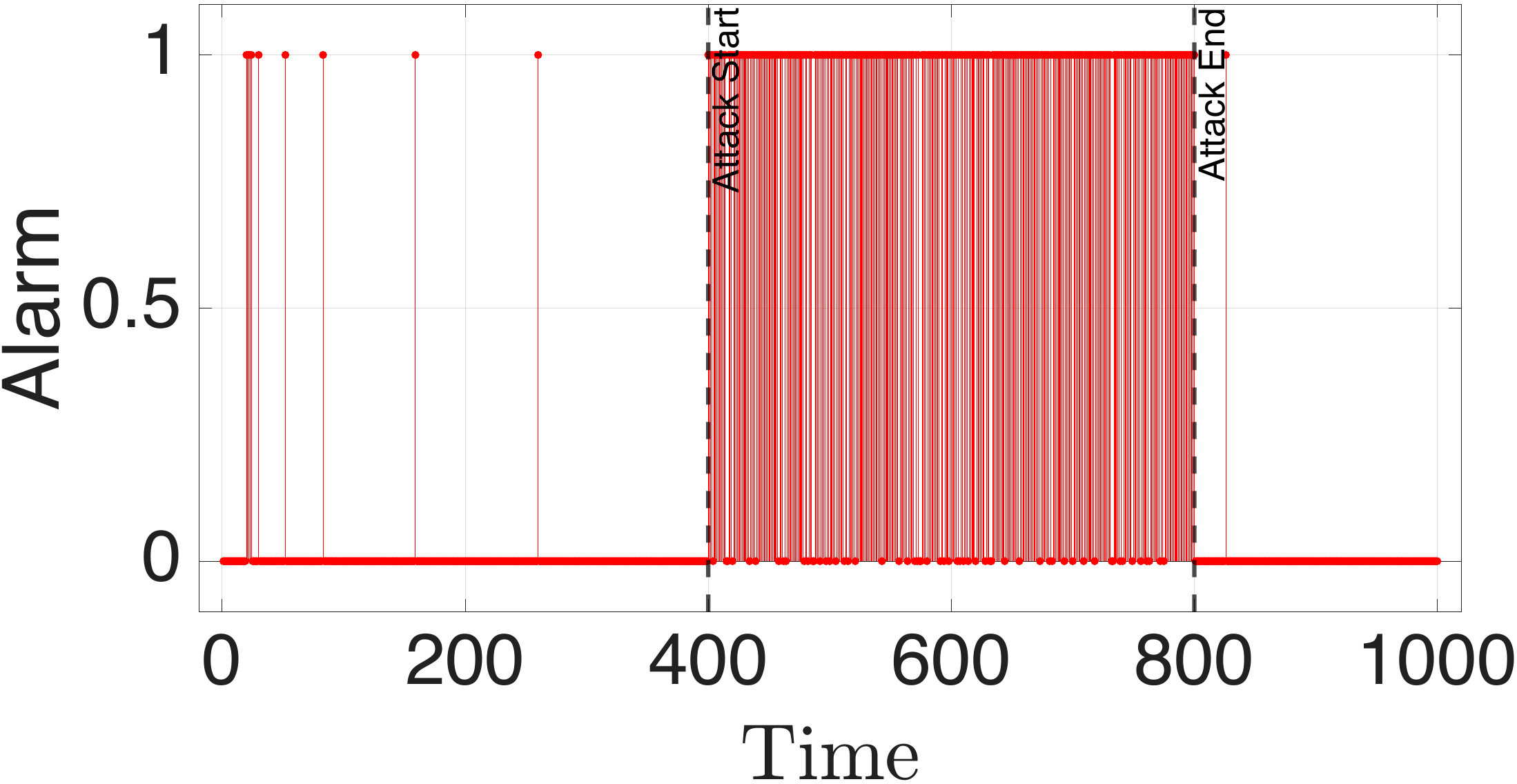}
    \end{subfigure}
    \caption{AD-CPS detector.}
    \label{subfig:AD-CPS}
\end{subfigure}

\begin{subfigure}{\linewidth}
    \centering
    \begin{subfigure}{0.48\linewidth}
        \includegraphics[width=\linewidth,height = 2.7cm]{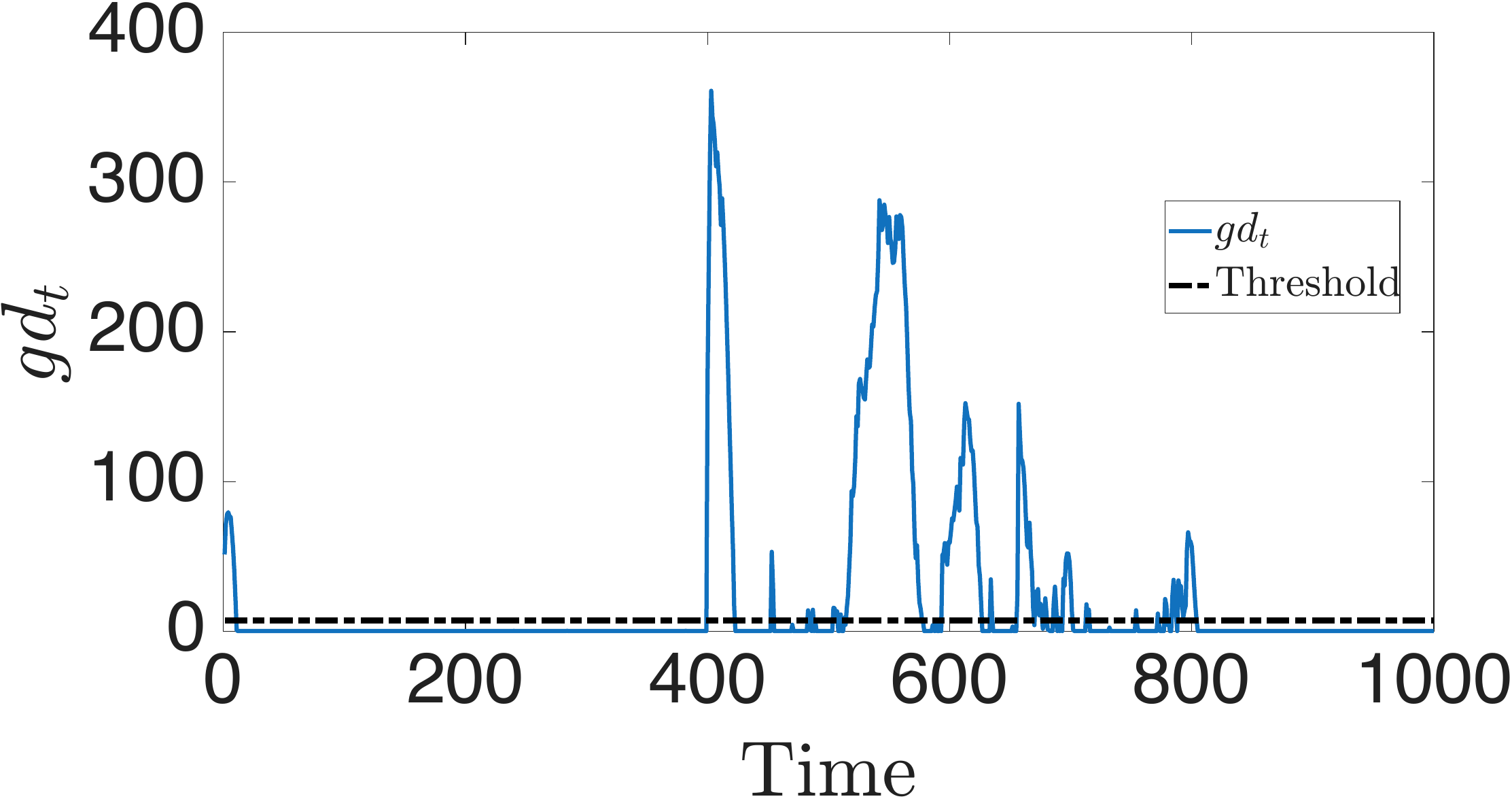}
    \end{subfigure}
    \hfill
    \begin{subfigure}{0.48\linewidth}
        \includegraphics[width=\linewidth,height = 2.56cm]{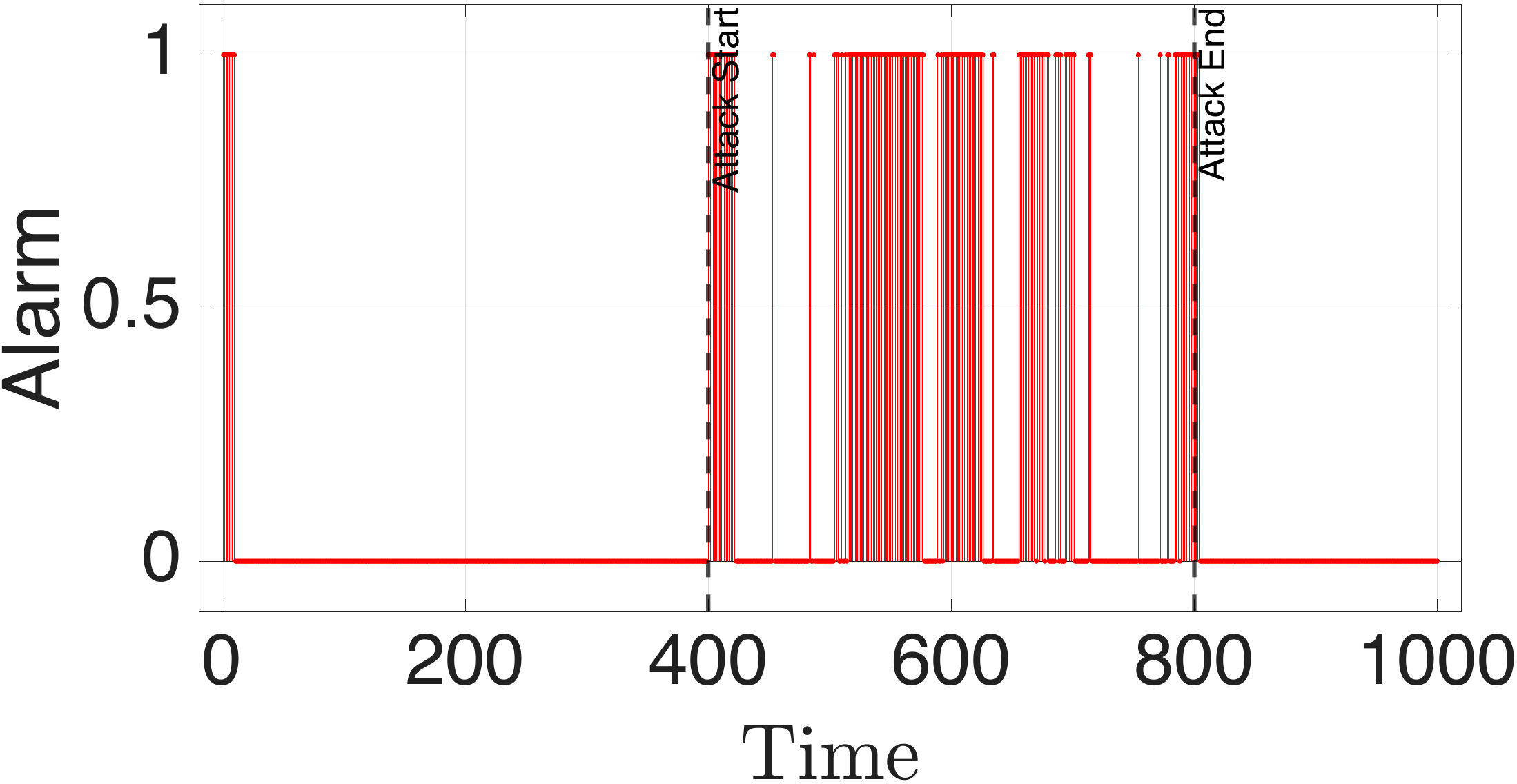}
    \end{subfigure}
    \caption{Optimal CUSUM-based detector.}
    \label{subfig:optimal CUSUM}
\end{subfigure}

\caption{Test signal and detection alarm using two different detection schemes. The simulations were conducted with \(\vbound=0.001\), \(\sbound_{\ol{\pnoise}}=0.01\), \(\sbound_a = 0.1\), and \(\window=20\).}
\label{fig:comparison}
\end{figure}

\textbf{Comparison with optimal CUSUM-based detector \cite{ref:AN-AT-AA-SD-23}. }We compared AD-CPS with an optimal watermarking-based CUSUM test. We conducted three experimental trials (with different seed values), and have shown the simulation corresponding to the best trial run in Figure \ref{fig:comparison}. Table \ref{tab:comparison} summarizes the comparison results for all the trial runs.

\begin{table}[htbp]
\centering

\caption{}
\label{tab:comparison}

\begin{adjustbox}{max width=0.9\linewidth}
\begin{tblr}{
  colspec = {l |cccccc},
  row{1,2} = {azure9},
  hline{1,3,Z} = {2pt},
  hline{2} = {1pt},
  column{2-7} = {c},
}
\SetTblrInner{rowsep=1pt}

& \SetCell[c=2]{c} Trail run 1 & &
  \SetCell[c=2]{c} Trail run 2 & &
  \SetCell[c=2]{c} Trail run 3 & \\

Algorithms &
\(\widehat{FPE}_{\window,\delta}\) &
\(\widehat{DR}_{\window,\delta}\) &
\(\widehat{FPE}_{\window,\delta}\) &
\(\widehat{DR}_{\window,\delta}\) &
\(\widehat{FPE}_{\window,\delta}\) &
\(\widehat{DR}_{\window,\delta}\)
\\

AD-CPS &
0.02 & 0.8454 & 0.0184 & 0.8554 & 0.01 & 0.8204\\

CUSUM-based detector \cite{ref:AN-AT-AA-SD-23} &
0.025063 & 0.2968 & 0.0256063 & 0.4788 & 0.027569 & 0.3915 \\

\end{tblr}
\end{adjustbox}

\end{table}
We highlight several key observations. The CUSUM-based test is specifically designed to achieve the quickest possible detection of attacks and is therefore optimal in that regard. However, the accuracy of AD-CPS, measured in terms of \(\wh{DR}_{\window,\delta}\) surpasses that of the CUSUM-based test. Moreover, the performance of the CUSUM-based test may deteriorate further relative to AD-CPS when the Gaussian noise assumption is violated.





\subsection{Computation of \(\log\Bigl(\frac{\epsilon^2}{k\vbound \ol{M}}\Bigr)\) in Theorem \ref{th:t2_error_incurred}} \label{subsec:fne compute}
\begin{figure}[htpb!]
\centering
 \includegraphics[scale=0.10]{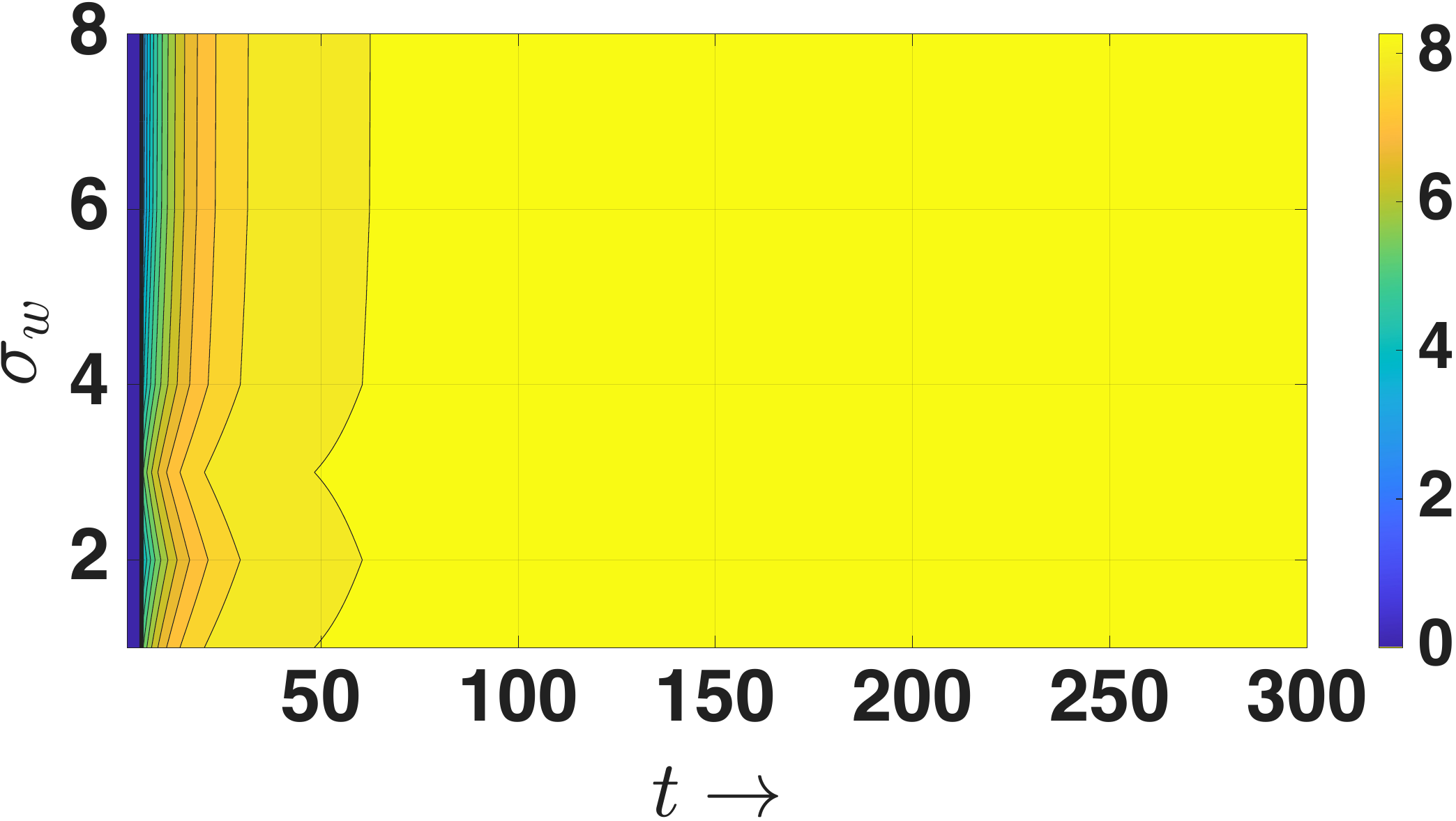}
 \caption{Plot of \(\log\Bigl(\epsilon^2/k\vbound \ol{M}\Big)\).}
 \label{fig:err_fne}
 \end{figure}

 In Figure \ref{fig:err_fne} we plotted \(\log\Bigl(\frac{\epsilon^2}{k\vbound \ol{M}}\Bigr)\) as a function of the process noise parameter \(\vbound\) and time \(t\). 
 We observe that for a fixed \(\vbound\), as time increases the value admitted by the quantity \(\log\Bigl(\frac{\epsilon^2}{k\vbound \ol{M}}\Bigr)\) increases, demonstrating that the false negative error (in Theorem \ref{th:t2_error_incurred}) decreases and saturates at \(2 \delta\). We now fix \(t\) and observed that as \(\vbound\) increases, the false negative error decreases, validating the assertion in Theorem \ref{th:t2_error_incurred}. 
\section{Conclusion}
This article establishes a finite-time algorithm (AD-CPS) with non-asymptotic guarantees. A test with a suitable threshold was designed to detect these attackers with high probability.  We investigated various aspects of AD-CPS through numerical experiments. In contrast to the setting studied in this article, we extended the setting to a partially observed case, thereby confirming that the algorithm can be easily generalized to a more general setting. Our detector is compared with an optimal watermarking CUSUM-based detector \cite{ref:AN-AT-AA-SD-23}, and we report comparable performance.

Future work may be directed towards improving the threshold used for the test so that stealthy attacks with small magnitudes can also be detected. It is worthwhile to relax the sub-Gaussian assumption on the process noise and design tests compatible with heavy-tail noise. Future direction may also include integrating \emph{quickest detection} in conjunction to the present setting.

\appendix
\section{Analysis of AD-CPS }
\label{appen:proofs}
We start analyzing AD-CPS (Algorithm \ref{alg:FTAAD}) by stating a concentration inequality that establishes an upper bound on the tails of sub-Gaussian random vectors. 
Recall that \(\st\) is a \(\Rbb^d\mbox{-}\)valued sub-Gaussian random vector if there exists a \(\sigma> 0\) such that for all \(\alpha \in \Rbb^d\),
\begin{equation}\label{eq:sG}
    \EE \expecof[\big]{\exp{(\alpha^{\top}\st)}} \leq \exp{\biggl(\frac{\sigma^2 \norm{\alpha}^2 }{2} \biggr)},
\end{equation}
and we use the notation \(x \sim \subG_d(\sbound)\) (with slight abuse of notation). 
The following result is useful:
\begin{prop}{\cite[Proposition 19.1]{ref:AS-BN-23}}\label{prop:au_con_in}
    Let \(\pnoise\) be a \(d\text{-}\)dimensional sub-Gaussian random vector with zero mean and covariance \(\var.\) Then there exists \(C>80\) such that for sufficiently small \(\delta \in ]0,1[\) we have
    \begin{align*}
        \PP \probof[\Big]{\abs{\norm{\pnoise} - \EE \expecof[\big]{\norm{\pnoise}}} \geq C \sqrt{\sbound_{\max}(\var) \ln{(1/\delta)}}} \leq \delta,
    \end{align*}
    where \(\sbound_{\max}(\var)\) is the maximum eigenvalue of \(\var\).
\end{prop}
\begin{lem}
    \label{lem:concen_ineq}
    Let \(\pnoise\) be an \(\Rbb^d\mbox{-}\)valued sub-Gaussian random vector with zero mean and covariance \(\var\). Let \(\sigma>0\) be such that \(\opnorm{\var} \leqslant \sigma\). Then there exists a \(k>0\) such that for sufficiently small \(\delta \in ]0,1[\), \(\PP \probof[\Big]{\norm{\pnoise} > \sqrt{k\sigma d\ln{(1/\delta))}}} \leqslant \delta.\)
\end{lem}
\begin{pf}
    Let us first compute \(\EE \expecof[\big]{\norm{\pnoise}}\). Note that 
    \begin{align*}
        \EE \expecof[\big]{\norm{\pnoise}}  \leq \sqrt{\EE \expecof[\Bigg]{\sum_{i=1}^d \pnoise_i^2}} 
         &= \sqrt{ \sum_{i=1}^d \EE \expecof[\big]{\pnoise_i^2} }\\ & = \sqrt{ \trace(\var)}
          \leq \sqrt{\sbound_{\max}(\var) d }.
    \end{align*}
    Moreover, we have 
    \(
        \abs{\norm{\pnoise} - \EE \expecof[\big]{\norm{\pnoise}}} \geq \norm{\pnoise} - \EE \expecof[\big]{\norm{\pnoise}} \geq \norm{\pnoise} - \sqrt{\sbound_{\max}(\var)d },
    \)
    and consequently, invoking Proposition \ref{prop:au_con_in}
    \begin{align*}
        &\PP \probof[\Big]{\norm{\pnoise} \geq \sqrt{\sbound_{\max}(\var)d } + C\sqrt{\sbound_{\max}(\var) \ln{(1/\delta)} }}\\
        & \leq \PP \probof[\Big]{\abs{\norm{\pnoise} - \EE \expecof[\big]{\norm{\pnoise}}} \geq C \sqrt{\sbound_{\max}(\var) \ln{(1/\delta)}}} \leq \delta.
    \end{align*}
     Pick \(k \in \Rbb\) such that \(k = C^2+1\). Since, \(\sbound_{\max}(\var) = \sbound\) and for sufficient small \(\delta \in ]0,1[\),
     \begin{align}\label{eq:ax_i}
        \sqrt{\sbound d } + C \sqrt{\sbound_{\max}(\var) \ln{(1/\delta)} } \leq \sqrt{k\sigma d \ln{(1/\delta)} },
     \end{align}
     which implies that
    \(
        \PP \probof[\Big]{\norm{\pnoise} > \sqrt{k\sigma d \ln{(1/\delta)}}} \leq \PP \probof[\bigg]{\norm{\pnoise} \geq \sqrt{\sbound d } + C\sqrt{\sbound \ln{\frac{1}{\delta}} }} \leq \delta.
    \)
    This proves the assertion.
\end{pf}
\begin{rem}
    \label{rem:sg_prop}
    The concentration inequality obtained in Lemma \ref{lem:concen_ineq} is tight if \(\pnoise\) is assumed to be \(d\text{-}\)dimensional Gaussian or sub-Gaussian with independent coordinates See \cite[Theorem 6.3.2]{ref:RV-MR18} for an elaborate discussion. 
\end{rem}
An immediate consequence of Lemma \ref{lem:concen_ineq} is that for any \(d\text{-}\)dimensional sub-Gaussian random vector \(\pnoise\) and small enough \(\delta \in \loro{0}{1}\), the set \(\aset[\big]{\norm{\pnoise} \leq \sqrt{k\sbound d \ln(1/\delta)}}\) is a high-probability set. If not specified, it is understood that all computations throughout the article are performed \emph{on} such high probability sets.
Let us now recall that the closed-loop CPS \eqref{eq:closed_lsys} and Doob decomposition of \(\st_t\) at time \(t\) in \eqref{eq:decomposition}:
\begin{align*}
    \stt{}{t} &= \mart_t + \predic_t + \stt{}{t-2}\nn\\
    & \stackrel{\mathclap{(\dagger)}}{=} \sum_{i=t-1}^t \pnoise_i + \sum_{i=t-1}^t \bigl( \Acl \stt{}{i-1} - \stt{}{i-1
    }\bigr) + \stt{}{t-2} \quad \PP\mbox{-}\text{almost surely}.
\end{align*}

Recall the definition of \(\sbound_0\) from \eqref{eq:iterations}. From \eqref{eq:decomposition}, note that  \(\mart_t = \sum_{i=t-1}^t \pnoise_i\) is a zero mean sub-Gaussian random vector on \(\Rbb^d\) with covariance given by the expression 
\(\var_{\mart_t} = \sum_{i=t-1}^t \var_{\pnoise}=2 \var_{\pnoise}\)
in view of independence of \((\pnoise_i)_{i=1}^t\),  and \(\opnorm{\var_{\mart_t}} \leq \sum_{i=t-1}^t \norm{\var_{\pnoise}} \leqslant 2 \sbound_{\pnoise}\) (follows from Assumption \ref{assum:sub_gaussian_noise}). 
Note that \((\mart_{t-1}, \mart_t)\) are \(d\mbox{-}\)dimensional sub-Gaussian random vectors with \(\opnorm{\var_{\mart_t}} \leq 2 \vbound\). Similarly, for a fixed \(t \in \Nz\) and \(s =t-1,t\), we see that \(\Acl^{t-s}\mart_s\) is a sub-Gaussian random vector on \(\Rbb^d\) with zero mean, and 
\begin{equation}
    \label{eq:aux_var}
    \Var\bigl(\Acl^{t-s}\mart_s\bigr) = \Acl^{t-s} \var_{\mart_s} \bigr(\Acl^{t-s} \bigr)^{\top}, 
\end{equation}
with \(\norm{\Var\bigl(\Acl^{t-s}\mart_s\bigr)} \leq 2 \vbound \norm{\Acl^{t-s}}\).
We have the following auxiliary lemmas:
\begin{lem}
    \label{lem:sub_Gauss}
    Consider the control system \eqref{eq:lsys1}/\eqref{eq:closed_lsys} along with its associated data. Suppose that Assumption \ref{assum:sub_gaussian_noise} holds, and consider Algorithm \ref{alg:FTAAD}. For each \(t \in \Nz\), define \(d_{t} \Let \frac{1}{2}\sum_{s =t-1}^t  \Acl^{t-s}\mart_s.\)
    Then \((d_{t})_{t \in \Nz}\) is a sequence of \(\Rbb^d\mbox{-}\)valued zero mean sub-Gaussian random vectors with
    \(
        \norm{\Var(d_{t})}  \leqslant \frac{\sbound_{d_t}}{4} 
    \)
    where \(\sbound_{d_t}\) and \(M\) are defined in \eqref{eq:aux_v1}.
\end{lem}
\begin{pf}
Observe that for \(s = t-1,t\), the norm of the covariance of \(d_{t}\) at time \(t\) is given by 
    \begin{align*}
        \norm{\Var (d_{t})} &= \left\lVert \frac{1}{4}\sum_{s =t-1}^t  \Var\bigl(\Acl^{t-s}\mart_s\bigr) + \frac{1}{4} \Cov\Bigl(\Acl\mart_{t-1}, \mart_t \Bigr) \right. \\ & \hspace{3cm}\left.+ \frac{1}{4} \Cov\Bigl(\mart_t,\Acl\mart_{t-1} \Bigr)\right\rVert.
    \end{align*}
    Note that \(\EE \expecof[\big]{\mart_s} = 0\). It is straightforward to compute the second and the third terms in the above equality, and they are given by
    \(\Cov\Bigl(\Acl\mart_{t-1}, \mart_t \Bigr) = \Acl \var_{\pnoise}\) and \(\Cov\Bigl(\mart_t,\Acl\mart_{t-1} \Bigr) = \var_{\pnoise} (\Acl){\top}\). Similarly, \(\sum_{s =t-1}^t  \Var\bigl(\Acl^{t-s}\mart_s\bigr) = \Acl \var_{\pnoise} (\Acl){\top} + 2 \var_{\pnoise}\). Define the quantities \(M \Let \Bigl(2 + 2 \norm{\Acl} + \norm{\Acl}^2 \Bigr)\) and \(\sbound_{d_t} \Let \vbound M\). Combining all the pieces, we obtain 
    \begin{align}\label{eq:var_sum}
        \norm{\Var (d_{t})} &\leq \frac{1}{4} \norm{\Acl \var_{\pnoise} (\Acl){\top} + 2 \var_{\pnoise} + \Acl \var_{\pnoise} + \var_{\pnoise} (\Acl){\top} }\nn\\
        & \leq \frac{ \vbound}{4} \Bigl(2 + 2 \norm{\Acl} + \norm{\Acl}^2 \Bigr)  = \frac{\sbound_{d_t}}{4}.
    \end{align}
    The proof is now complete.
\end{pf}
%
\begin{lem}
    \label{lem:a_sub_Gauss}
    Suppose that the hypotheses of Lemma \ref{lem:sub_Gauss} hold. Define \(\ol{d}_{t} \Let d_{t} - \mart_t.\)
    Then \((\ol{d}_{t})_{t \in \Nz}\) is a sequence of \(\Rbb^d\mbox{-}\)valued zero mean sub-Gaussian random vectors with \(\opnorm{\Var(\ol{d}_{t})} \leqslant  \ol{\sbound}\) where \(\ol{\sbound} \) and \(\ol{M}\) are defined in \eqref{eq:aux_v1}.
\end{lem}
\begin{pf}
    The proof proceeds along the lines of Lemma \ref{lem:sub_Gauss}. Observe that \(\norm{\Var(\ol{d}_{j,t})} \leq \norm{\Var(d_{j,t})} + \norm{\var_{\mart_t}} + \norm{\Cov(d_{j,t},\mart_t))} + \norm{\Cov(\mart,d_{j,t})} \stackrel{\mathclap{(\star)}}{\leq} \biggl(2 + \frac{M}{4}\biggr) \vbound +\norm{\Cov(d_{t},\mart_t)} + \norm{\Cov(\mart_t,d_{t})}\),
    where \((\star)\) follows from Lemma \ref{lem:sub_Gauss}. Let us focus our attention on the second and third terms in \((\star)\). From Assumption \ref{assum:sub_gaussian_noise} we have \(\EE \expecof[\big]{\ol{d}_{t}} = 0\). Observe that \(\Cov(d_{t},\mart_t) = \frac{1}{2}\Acl \var_{\pnoise}+ \var_{\pnoise}\) and similarly, \(\Cov(\mart_t,d_{t}) = \frac{1}{2}\var_{\pnoise}(\Acl){\top} + \var_{\pnoise}\). This implies that \(\norm{\Cov(d_{j,t},\mart_t)} = \norm{\Cov(\mart_t,d_{t})} \leq \biggl( 1+ \frac{\norm{\Acl}}{2} \biggr)\vbound \).
    Combining everything, we get \(\norm{\Var(\ol{d}_{j,t})} \leq \biggl(2 + \frac{M}{4}\biggr) \vbound + \biggl( 1+ \frac{\norm{\Acl}}{2} \biggr)2 \vbound = \vbound \biggl(4 + \frac{M}{4} + \norm{\Acl} \biggr) = \vbound \ol{M} \teL \ol{\sbound} .\)
    where \(\ol{M} \Let \bigl(4 + \frac{M}{4} + \norm{\Acl} \bigr)\).
    This completes the proof.
\end{pf}
Recall the set \eqref{eq:good_set}, 
\begin{equation*}
    \goodset_t \Let  \aset[\Bigg]{\Bigg \lVert\frac{1}{2}\sum_{s =t-1}^t  \Acl^{t-s}\mart_s -\mart_t \Bigg \rVert \leqslant \sqrt{k \ol{\sbound} d \ln{\bigl(1/\delta \bigr)}}}. 
\end{equation*}
In view of Lemma \ref{lem:a_sub_Gauss} we observe that for each \(t \in \Nz\), \(\goodset_t\) is a ``high probability set''.
\begin{lem}
\label{lem:SG_good_set}
Suppose that the hypotheses of Lemma \ref{lem:sub_Gauss} hold, and consider Algorithm \ref{alg:FTAAD}. Then \(\PP\probof[\big]{\goodset_t} \geqslant 1- \delta.\)
\end{lem}
\begin{pf}
    From Lemma \ref{lem:a_sub_Gauss}, \(\ol{d}_t = \frac{1}{2}\sum_{s =t-1}^t  \Acl^{t-s}\mart_s -\mart_t \) is a \(d\)-dimensional sub-Gaussian random vector with parameter \(\ol{\sbound}\). 
    The proof is straightforward, invoking Lemma \ref{lem:concen_ineq}. 
\end{pf}
\section{Proof of Theorem \ref{th:error_incurred}}\label{appen:t1_error}
\begin{lem}
    \label{lem:predict_upperBnd}
    Suppose that the hypotheses of Lemma \ref{lem:sub_Gauss} hold. Define the quantity \(q_{t} \Let \frac{1}{2}\sum_{s = t-1}^t\Acl^{t-s}\predic_s - \predic_t.\)
    Then we have the following assertions:
    \begin{enumerate}[label=\textup{(\ref{lem:predict_upperBnd}-\alph*)}, leftmargin=*, widest=b, align=left]
        \item \label{it:sub_Gauss_a}\(\PP \bigl(\norm{\st_t}  \leq \sqrt{k \sbound_{\st_t}d \ln(1/\delta) }\bigr) \geq 1 - \delta \), where \(\sbound_{\st_{t}} = \sbound_0 \norm{\Acl^{t}}^2 + \vbound \sum_{j=0}^{t-1} \norm{\Acl^{t-1-j}}^2\);
        
        \item \label{it:lisan_al_gaib} \(\norm{q_{t}} \leq \norm{\Acl - \identity{d}}\Bigl(\sum_{s =t-1}^t \norm{\Acl^{t-s}}+1\Bigr)\sum_{i=t-1}^t \norm{\st_{i-1}}.\)
    \end{enumerate}
\end{lem}
\begin{pf}
    The first part of the proof directly follows by invoking Lemma \ref{lem:concen_ineq} and from the fact that \(\st_t\) for each \(t \in \Nz\), is \(d\text{-}\)dimensional sub-Gaussian random vector.
    
    For the second part of the proof, a direct application of triangle inequality yields 
    \begin{align*}
        &\norm{q_{t}} \leq \sum_{s =t-1}^t \norm{\Acl^{t-s}}\sum_{i=t-1}^s \norm{\Acl - \identity{d}}\norm{\st_{i-1}} + \sum_{k=t-1}^t \norm{\Acl^{t-s}}\norm{\st_{k-1}}\\
        &\leq   \norm{\Acl - \identity{d}}\biggl(\sum_{s =t-1}^t \norm{\Acl^{t-s}}+1\biggr)\sum_{i=t-1}^t  \norm{\st_{i-1}}.
    \end{align*}
    Define the quantity \(\ol{h} \Let \norm{\Acl - \identity{d}}\Bigl(\sum_{s =t-1}^t \norm{\Acl^{t-s}}+1\Bigr) \). Then we get
    \(\norm{q_{t}} \leq \ol{h}\sum_{i=t-1}^t  \norm{\st_{i-1}}\). This completes the proof.
\end{pf}
Now we give the proof of Theorem \ref{th:error_incurred} in full detail. 
\begin{pf}[Proof of Theorem \ref{th:error_incurred}]
Let us recall
\begin{align*}
\begin{cases}
    q_{t} \Let \frac{1}{2}\sum_{s =t-1}^t\Acl^{t-s}\predic_s - \predic_t,\,
    \ol{d}_{t} = \frac{1}{2}\sum_{s =t-1}^t\Acl^{t-s}\mart_s - \mart_t,\\
    \adver_t =  \frac{1}{2}\sum_{s =t-1}^t\Acl^{t-s}\donoise{s} - \donoise{t} -\frac{1}{2}\sum_{s=t-1}^t (\Acl^{t-s} - \identity{d})\donoise{t-2},
\end{cases}
\end{align*}
and the \emph{test-signal} from \eqref{eq:real} by
\begin{align*}
    \tstat_t &\Let \frac{1}{j}\sum_{s \in \indexset_j}\Acl^{t-s}\sm{}{s} - \sm{}{t}- \frac{1}{2}\sum_{s=t-1}^t (\Acl^{t-s} - \identity{d})\sm{}{t-2} \nn\\
    & = \ol{d}_t + q_t + \adver_t
\end{align*}
Fix \(t \in \Nz\); define the set 
\(
    \thonest_t \Let \aset[\big]{\sm{}{t} \text{ has \(2\)-step honest memory}},
\)

and observe that on sets \(\goodset_t\) and \(\thonest_t\), the conditional probability of interest is given by
\begin{align*}
    &\PP \cprobof[\Big]{\aset[\big]{\text{Algorithm \ref{alg:FTAAD} declares }\sm{}{t} \text{ as corrupt}} \cap \thonest_t \cap \goodset_t\given \hypothesis_{0,t} }\\
    & = \PP \cprobof[\Big]{\aset[\big]{\norm{\tstat_t} >\athr}\cap \aset[]{\donoise{t-2}=\donoise{t-1}=0}  \cap \goodset_t \given \hypothesis_{0,t}}
\end{align*}
Observe that given \(\hypothesis_{0,t}\),
\begin{align*}
    &\aset[\big]{\norm{\tstat_t} >\athr}\cap \aset[]{\donoise{t-2}=\donoise{t-1}=0}  \cap \goodset_t = \aset[\Bigg]{\norm{\ol{d}_t + q_t } > \athr}\cap \goodset_t,
\end{align*}
implying that \(\PP \cprobof[\Big]{\aset[\big]{\norm{\tstat_t} >\athr}\cap \aset[]{\donoise{t-2}=\donoise{t-1}=0}  \cap \goodset_t \given \hypothesis_{0,t}}\) can be estimated as
\begin{align*}
    &\PP \cprobof[\Bigg]{\aset[\big]{\norm{\tstat_t} >\athr}\cap \aset[]{\donoise{t-2}=\donoise{t-1}=0}  \cap \goodset_t \given \hypothesis_{0,t}}\\
    & = \PP \probof[\Bigg]{\aset[\Bigg]{\norm{\ol{d}_t + q_t} > \athr}\cap \goodset_t} \\
    & = \PP \probof[\Bigg]{\aset[\Bigg]{\norm{\ol{d}_t + q_t} > \Bigl(\sqrt{2} + \sqrt{\ol{M}}\Bigr) \sqrt{k \vbound d \ln(1 / \delta)}}\cap \goodset_t} \\
    & \leq \PP \probof[\Bigg]{\aset[\Bigg]{\norm{\ol{d}_t}  >  \sqrt{k \vbound \ol{M} d \ln(1 / \delta)}}\cap \goodset_t} \\& \hspace{2cm} + \PP \probof[\Bigg]{\aset[\Bigg]{\norm{q_t}  >  \sqrt{2 k \vbound  d \ln(1 / \delta)}}\cap \goodset_t}\\
    & \leq \delta + \PP \probof[\Big]{\norm{q_t}  >  \sqrt{2 k \vbound  d \ln(1 / \delta)}} \quad \text{by applying Lemma \ref{lem:a_sub_Gauss}},\\
    & \leq \delta + \PP \probof[\Bigg]{\ol{h}\sum_{i=t-1}^t  \norm{\st_{i-1}}  >  \sqrt{2 k \vbound  d \ln(1 / \delta)}}\;\text{from Lemma \ref{it:lisan_al_gaib}},
    \\
    & \leq \delta + \sum_{i=t-1}^t \PP \probof[\Bigg]{\norm{\st_{i-1}}  >  \frac{\sqrt{2 k \vbound  d \ln(1 / \delta)}}{2 \ol{h}}}\\
    & \leq \delta + \sum_{i=t-1}^t \frac{2 \ol{h}}{\sqrt{2 k \vbound  d \ln(1 / \delta)}} \EE \expecof[\big]{\norm{\st_{i-1}}},
\end{align*}
where the last inequality follows from Markov inequality \cite[\S5, pg. 46]{ref:Shi-16}.
Let us now focus on the term \(\EE \expecof[\big]{\norm{\st_{i-1}}}\) where \(i=t-1,t\). We first observe that \((\st_t)_{t \in \Nz}\) is  Markov process. From the Markovian property and the `Telescopic property' of conditional expectation, 
\(
    \EE \expecof[\big]{\norm{\st_{i-1}}} = \EE \expecof[\Big]{ \EE \cexpecof[\big]{\norm{\st_{i-1}} \given \st_0} } 
\), where the inner conditional expectation \(\EE \cexpecof[\big]{\norm{\st_{i-1}} \given \st_0}\) can be expanded as
\begin{align*}
    &\EE \cexpecof[\big]{\norm{\st_{i-1}} \given \st_0} = \EE \left[ \EE \cexpecof[\big]{\norm{\st_{i-1}} \indic{K}(\st_{t-2}) \given \st_{i-2}}\right. \\& \hspace{1cm} + \left. \EE \cexpecof[\big]{\norm{\st_{i-1}} \indic{K^{\complift}}(\st_{t-2}) \given \st_{i-2}} \Big|  \st_0 \right]   \\
    & \stackrel{\mathclap{(\star)}}{\leq} \EE \cexpecof[\Big]{\indic{K}(\st_{i-2}) \EE \cexpecof[\big]{\norm{\st_{i-1}}  \given \st_{i-2}} + \gamma \norm{\st_{i-2}} \given \st_0}, \\
    & \leq \beta \PP \cprobof[\big]{\st_{i-2} \in K \given \st_0} + \gamma \EE \cexpecof[\big]{ \norm{\st_{i-2}} \given \st_0},
\end{align*}
where \((\star)\) follows from hypotheses of Theorem \ref{th:error_incurred}.
Proceeding iteratively, we obtain
\(
     \EE \cexpecof[\big]{\norm{\st_{i-1}} \given \st_0} = \gamma^{i-1} \norm{\st_0} + \beta \sum_{j=0}^{i-2} \gamma^{i-1-j}
\), and 
\begin{align*}
    \EE \expecof[\big]{\norm{\st_{i-1}}} &= \EE \expecof[\Big]{ \EE \cexpecof[\big]{\norm{\st_{i-1}} \given \st_0} } 
    & = \gamma^{i-1} \norm{\st_0} + \beta \sum_{j=0}^{i-2} \gamma^{i-1-j}.
\end{align*}
Putting all the pieces together, we write
\begin{align*}
    &\PP \cprobof[\Bigg]{\aset[\big]{\norm{\tstat_t} >\athr}\cap \aset[]{\donoise{t-2}=\donoise{t-1}=0}  \cap \goodset_t \given \hypothesis_{0,t}}\\
    & \leq \delta + \sum_{i=t-1}^t \frac{2 \ol{h}}{\sqrt{2 k \vbound  d \ln(1 / \delta)}} \EE \expecof[\big]{\norm{\st_{i-1}}} \\
    & \leq \delta + \sum_{i=t-1}^t \frac{2 \ol{h} \Bigl(\gamma^{i-1} \norm{\st_0} + \beta \sum_{j=0}^{i-2} \gamma^{i-1-j}\Bigr)}{\sqrt{2 k \vbound  d \ln(1 / \delta)}}.
\end{align*}
This completes the proof.
\end{pf}
\section{Proof of Theorem \ref{th:t2_error_incurred}}\label{appen:t2_error}
The following lemma is useful in the proof of Theorem \ref{th:t2_error_incurred}.
\begin{lem}
    \label{lem:a_result}
    Suppose that the hypotheses of Lemma \ref{lem:sub_Gauss} hold. Define the quantities \(\ol{b}_{t} \Let \mart_t - \frac{1}{2}\sum_{s=t-1}^t\Acl^{t-s}\mart_s \) and \(\phi \Let \vbound \ol{M}\).
    Then for each \(t \in \Nz\), we have
    \(
        \norm{\Var\bigl(\ol{b}_{t}\bigr)} \leq \phi
    \).
\end{lem}
%
The proof is straightforward, so we omit it.
Now we proceed to the proof of Theorem \ref{th:t2_error_incurred}. We clinically adopt the notation used in Lemma \ref{lem:predict_upperBnd} and Lemma \ref{lem:a_result} above.
\begin{pf}[Proof of Theorem \ref{th:t2_error_incurred}]
Recall from \eqref{eq:real}, the definition of \(\adver_t\). Fix \(t \in \Nz\). On the event \(\aset[\big]{\norm{\adver_t}\geq \thres_t}\) and given \(\hypothesis_{1,t}\) we are interested in the probability
\begin{align*}
    \PP \cprobof[\Big]{\aset[\big]{\norm{\T_t \leq \athr}} \cap \aset[\big]{\norm{\adver_t} \geq \thres_t} \given \hypothesis_{1,t}}.
\end{align*}
A careful scrutiny reveals 
\(
    \norm{\T_t} \geq \norm{\adver_t} - \norm{q_{t}} - \norm{\ol{b}_{t} - s_{t}}
\).
This implies that 
\begin{align*}
    &\PP \cprobof[\Big]{\aset[\big]{\norm{\T_t \leq \athr}} \cap \aset[\big]{\norm{\adver_t} \geq \thres_t} \given \hypothesis_{1,t}}\\
    & \leq \PP \cprobof[\Big]{ \aset[\Big]{\norm{\adver_t} - \norm{q_{t}} - \norm{\ol{b}_{t} - s_{t}} \leq \athr} \cap \aset[\big]{\norm{\adver_t} \geq \thres_t} \given \hypothesis_{1,t} }\\
    &= \PP \cprobof[\Big]{\aset[\Big]{\norm{\ol{b}_{t} - s_{t}}  \geq \norm{\adver_t} - \norm{q_{t}} -  \athr } \cap \aset[\big]{\norm{\adver_t} \geq \thres_t} \given \hypothesis_{1,t} }\\
    & \leq \PP \cprobof[\bigg]{\norm{\ol{b}_{t} - s_{t}} + \norm{q_{t}} \geq \thres_{t} -  \athr \given \hypothesis_{1,t}}.
\end{align*}
Recall from \eqref{eq:adv_lb}, the definition of \(\lbt_t(\Acl,d,\vbound,\delta)\). Let Assumption \ref{assum:ad_assum} be in force. Then we have \(\thres_{t} -  \athr = \lbt_t(\Acl,d,\vbound,\delta) + \athr - \athr = \lbt_t(\Acl,d,\vbound,\delta)\). Let \(\ol{h} \Let \norm{\Acl - \identity{d}}\Bigl(\sum_{s =t-1}^t \norm{\Acl^{t-s}}+1\Bigr) \). The following computations reveal that
\begin{align*}
    &\PP \cprobof[\bigg]{\norm{\ol{b}_{t} - s_{t}} + \norm{q_{t}} \geq \thres_{t} -  \athr \given \hypothesis_{1,t}}\\
    & \leq \PP \cprobof[\bigg]{\norm{\ol{b}_{t} - s_{t}} + \ol{h} \sum_{i=t-1}^t \norm{\st_{i-1}} \geq \lbt_t(\Acl,d,\vbound,\delta) \given \hypothesis_{1,t}}\\
    & = \PP \Biggl(\norm{\ol{b}_{t} - s_{t}} + \ol{h} \sum_{i=t-1}^t \norm{\st_{i-1}} \geq \sqrt{d \vbound \ol{M}} \cdots \\ & \hspace{1cm} \cdots + \ol{h}\sum_{i=t-1}^t \Bigl( \sqrt{d \sbound_{\st_{i-1}}} + \sqrt{k\sbound_{\st_{i-1}} \ln{(1/\delta)}} \Bigr) \Biggr)\\
    & \leq \PP \cprobof[\Bigg]{\norm{\ol{b}_{t} - s_{t}} \geq \sqrt{d \vbound \ol{M}}  + \frac{\ol{h}}{2}\sum_{i=t-1}^t \sqrt{k\sbound_{\st_{i-1}} \ln{(1/\delta)}}} \\ &  + \PP \Bigg( \ol{h}  \sum_{i=t-1}^t \norm{\st_{i-1}} \geq \ol{h}  \sum_{i=t-1}^t \sqrt{d\sbound_{\st_{i-1}}} \cdots \\& \hspace{1cm} \cdots + \frac{\ol{h}}{2}\sum_{i=t-1}^t \sqrt{k\sbound_{\st_{i-1}} \ln{(1/\delta)}}\Bigg)
\end{align*}
Denote by \(\eps \Let \frac{\ol{h}}{2}\sum_{i=t-1}^t \sqrt{k \sbound_{\st_{i-1}} \ln{(1/\delta)}}\) and \(\phi \Let \vbound \ol{M} \). The above inequality is then expanded as
\begin{align*}
    &\PP \cprobof[\Bigg]{\norm{\ol{b}_{t} - s_{t}} \geq \sqrt{d \phi}  + \eps} \\& \hspace{0.5cm} + \sum_{i=t-1}^t \PP \probof[\bigg]{ \norm{\st_{i-1}} \geq \sqrt{d\sbound_{\st_{i-1}}} + \frac{1}{2} \sqrt{k\sbound_{\st_{i-1}} \ln{(1/\delta)}}}\\
    & \leq \exp{\Biggl( - \frac{\eps^2}{k \phi}\Biggr)} + \sum_{i=t-1}^t \PP \Bigg( \norm{\st_{i-1}} \geq \sqrt{d\sbound_{\st_{i-1}}} \cdots \\&\hspace{2cm} \cdots + \frac{1}{2} \sqrt{k \sbound_{\st_{i-1}} \ln{(1/\delta)}}\Bigg) \; \text{from Lemma \ref{lem:a_result}}.
\end{align*}
Let us compute \(\PP \probof[\Big]{ \norm{\st_{i-1}} \geq \sqrt{d\sbound_{\st_{i-1}}} + \frac{1}{2} \sqrt{k \sbound_{\st_{i-1}} \ln{(1/\delta)}}}\). Recall that for each \(t \in \Nz\) and from Assumption \ref{assum:sub_gaussian_noise}, \(\st_t\) is \(d\)-dimensional sub-Gaussian random vector. Invoking Lemma \ref{lem:concen_ineq} and Lemma \ref{it:sub_Gauss_a}, for \(\ol{\epsilon}>0\) we get 
\(
    \PP \probof[\Big]{ \norm{\st_{i-1}} \geq \sqrt{d\sbound_{\st_{i-1}}} +\ol{\epsilon}}
     \leq \exp{\Bigl(- \frac{\ol{\epsilon}^2}{k\sbound_{\st_{i-1}}} \Bigr)}
\). Pick \(\ol{\epsilon} = \frac{1}{2} \sqrt{k \sbound_{\st_{i-1}} \ln{(1/\delta)}}\) and 
\(
    \PP \probof[\bigg]{ \norm{\st_{i-1}} \geq \sqrt{d\sbound_{\st_{i-1}}} + \frac{1}{2} \sqrt{k\sbound_{\st_{i-1}} \ln{(1/\delta)}}} \leq \delta. 
\)
Collecting the pieces together, we write 
\begin{align*}
    &\PP \cprobof[\bigg]{\norm{\ol{b}_{t} - s_{t}} + \norm{q_{t}} \geq \thres_{t} -  \athr \given \hypothesis_{1,t}}\\
    & \leq \exp{\Biggl( - \frac{\eps^2}{k \phi}\Biggr)} + \sum_{i=t-1}^t \delta
     = \exp{\Biggl( - \frac{\eps^2}{k \phi}\Biggr)} + 2\delta,
\end{align*}
implying that 
\(
    \PP \cprobof[\Big]{\aset[\big]{\norm{\T_t \leq \athr}} \cap \aset[\big]{\norm{\adver_t} \geq \thres_t} \given \hypothesis_{1,t}}
    \leq \exp{\biggl( - \frac{\eps^2}{k \phi}\biggr)} + 2\delta.
\)
This concludes the proof.
\end{pf}

\bibliographystyle{ieeetr}
\bibliography{ref}

\end{document}